\documentclass[11pt,a4paper]{amsart}
\usepackage{color}
\usepackage{amssymb}
\usepackage{accents}
\def\R{\mathbb{R}}
\def\m1{{I\!\!M}}
\newcommand{\CC}{\mathbb{C}}
\newcommand{\NN}{\mathbb{N}}
\newcommand{\RR}{\mathbb{R}}
\newcommand{\rdue}{\R^2}
\newcommand{\ci}{\mathbb{C}}

\renewcommand{\to}{\rightarrow}
\newcommand{\pa}{\partial}

\newcommand{\ino}{\int_{\Omega}}

\newcommand{\ii}{\infty}
\newcommand{\eps}{\varepsilon}
\newcommand{\dt}{\delta}

\newcommand{\al}{\alpha}

\newcommand{\sg}{\sigma}

\newcommand{\om}{\Omega}

\newcommand{\lm}{\lambda}
\newcommand{\ov}[1]{\overline{#1}}
\newcommand{\uv}[1]{\underline{#1}}

\newcommand{\ainf}{\mbox{as\;}\;n\to+\infty}

\newcommand{\fo}{\forall}
\newcommand{\e}[1]{{\,\dsp e^{\dsp #1}}}

\newtheorem{theorem}{Theorem}[section]
\newtheorem{proposition}[theorem]{Proposition}
\newtheorem{lemma}[theorem]{Lemma}
\newtheorem{corollary}[theorem]{Corollary}
\newtheorem{remark}[theorem]{Remark}
\newtheorem{definition}[theorem]{Definition}
\newcommand{\finedim}{\hspace{\fill}$\square$}

\newcommand{\fineproof}{\hspace{\fill}$\square$}

\newcommand{\rife}[1]{(\ref{#1})}

\newcommand{\sscp}{\scriptscriptstyle}
\newcommand{\dsp}{\displaystyle}

\newcommand{\brm}{\begin{remark}\rm}
\newcommand{\erm}{\end{remark}}
\newcommand{\bdf}{\begin{definition}\rm}
\newcommand{\edf}{\end{definition}}
\newcommand{\bte}{\begin{theorem}}
\newcommand{\ete}{\end{theorem}}
\newcommand{\bpr}{\begin{proposition}}
\newcommand{\epr}{\end{proposition}}
\newcommand{\ble}{\begin{lemma}}
\newcommand{\ele}{\end{lemma}}
\newcommand{\bco}{\begin{corollary}}
\newcommand{\eco}{\end{corollary}}
\newcommand{\beq}{\begin{equation}}
\newcommand{\eeq}{\end{equation}}
\newcommand{\bdm}{\begin{displaymath}}
\newcommand{\edm}{\end{displaymath}}

\newcommand{\graf}[1]{\left\{\begin{array}{ll}#1\end{array}\right.}

\def\sideremark#1{\ifvmode\leavevmode\fi\vadjust{\vbox to0pt{\vss% the remark
 \hbox to 0pt{\hskip\hsize\hskip1em%                          will appear only
 \vbox{\hsize2.1cm\tiny\raggedright\pretolerance10000%          on the side
  \noindent #1\hfill}\hss}\vbox to15pt{\vfil}\vss}}}%

\begin{document}
\numberwithin{equation}{section}
\parindent=0pt
\hfuzz=2pt
\frenchspacing

\title[A Singular Liouville-type equation and the Alexandrov isoperimetric inequality]{On a singular Liouville-type equation and the Alexandrov isoperimetric inequality.}

\author[D.B. \& D.C.]{Daniele Bartolucci$^{(1,\dag)}$, Daniele Castorina$^{(2,\ddag)}$}

\thanks{2010 \textit{Mathematics Subject classification:} 35B45, 35J75, 35R05, 35R45, 30F45, 53B20. }

\thanks{$^{(1)}$Daniele Bartolucci, Department of Mathematics, University
of Rome {\it "Tor Vergata"}, \\  Via della ricerca scientifica n.1, 00133 Roma,
Italy. e-mail:bartoluc@mat.uniroma2.it}

\thanks{$^{(2)}$Daniele Castorina, Dipartimento di Matematica, Universit\`a di Padova, \\
Via Trieste 63, 35121 Padova, Italy. e-mail: castorin@math.unipd.it}

\thanks{$^{(\dag)}$Research partially supported by FIRB project {\sl
Analysis and Beyond} and by PRIN project 2012, ERC PE1\_11,
"Variational and perturbative aspects in nonlinear differential problems"}

\thanks{$^{(\ddag)}$Research partially supported by project {\sl Bando Giovani Studiosi 2013 - Universit\`a di Padova - GRIC131695}}

\begin{abstract}
We obtain a generalized version of an inequality, first derived by C. Bandle in the analytic setting, for weak subsolutions of a singular Liouville-type equation.
As an application we obtain a new proof of the Alexandrov isoperimetric inequality on  singular abstract surfaces. 
Interestingly enough, motivated by this geometric problem, we obtain a seemingly new characterization of local metrics on 
Alexandrov's surfaces of bounded curvature. At least to our knowledge, the characterization of the equality case in the 
isoperimetric inequality in such a weak framework is new as well.
\end{abstract}

\maketitle
{\bf Keywords}: Singular Liouville-type equations, Alexandrov's Isoperimetric inequality, Surfaces of Bounded Curvature.
\section{Introduction}
Let $\om\subset \RR^2$ be an  open, smooth and bounded domain, $K$ a measurable function on $\om$, $\omega$ be a signed measure of
bounded total variation in ${\om}$ and $\omega=\omega_{+}-\omega_{-}$ be its
Jordan decomposition, that is, for a
Borel set $E\subseteq {\om}$, $\omega_\pm(E)=\sup\limits_{U\subset E}(\pm\omega(U))$.
Then $\omega_{\pm}$ are non negative and mutually orthogonal measures of bounded total variation on ${\om}$ and we define
$f=f_+-f_-$, where $f_+$ and $f_-$ are two superharmonic functions constructed as follows,
\beq\label{frep}
f_{\pm}(x)=h_{\pm}(x)+\ino G(x,y) d\omega_{\pm}(y),
\eeq
where $h_{\pm}$ are harmonic in $\om$.  Here $G(x,y)$ denotes the Green's function
of $-\Delta$ in $\om$. We are concerned with some quantitative
estimates for subsolutions of the Liouville-type equation,
\beq\label{1n}
-\Delta u =   2K e^{f}\e{u}\;\; \mbox{ in }\;\;\om.
\eeq

By assuming $\om$ simply connected, $\pa\om$ analytic, $u$ an analytic and $C^{0}(\ov{\om})$ subsolution of \eqref{1n} with
$K(x)\equiv K_0$ in $\om$ for some $K_0\geq 0$,
in a pioneering paper \cite{Band0}, C. Bandle proved that,
\beq\label{band}
L^2(\pa \om)\geq \left(4\pi - \omega_{+}(\om)-K_0M(\om)\right)M(\om),
\eeq
where,
$$
L(\pa \om)=\int\limits_{\pa \om} e^{\frac{f+u}{2}}d\ell,\quad\mbox{and}\quad
M( \om )=\int\limits_{ \om}e^{f+u}dx.
$$
Here and in the rest of this paper $d\ell$ and $dx$ will denote the $1$-dimensional and $2$-dimensional Hausdorff measures respectively.

The inequality \rife{band} is sharp and in \cite{Band0} the case where the equality holds is characterized as well, see also \cite{Band}.
Actually \rife{band} admits a beautiful geometric interpretation in terms of the Alexandrov's isoperimetric inequality \cite{Ale},
as discussed in \cite{Band0} and more extensively in \cite{Band}. As we will see later on,
the original geometric setting of the problem in terms of singular isothermal coordinates \cite{Res2}, suggests that \rife{band}
should hold in a more general form. This is our motivation and indeed our main aim is to obtain a
generalized version of \rife{band} in a weak framework. To state our result, we need some definitions first.

\bdf\label{simpn}{\it
We say that $E\subset \R^2$ is a \uv{simple} domain, if
it is an open and bounded domain whose boundary $\pa E$ is the support of
a rectifiable Jordan curve.
We will also say that $E\subset \R^2$ is a \uv{regular} domain if it is a connected,
open and bounded domain whose boundary $\pa E$ is the union of finitely many rectifiable
Jordan curves.}
\edf

\bdf\label{def12} {\it Let $S\subset \om$ be a finite set. We say that
$$
f\in L^{p,{\sscp \rm loc}}_{\rm loc}\left(\om\setminus S\right)\mbox{ or either }
u\in W^{2, p, {\sscp \rm loc}}_{\rm loc}\left(\om\setminus S\right),\;\mbox{ for some }p>2,
$$
if for each open and relatively compact set $U\Subset \om\setminus S$  there exists $p=p_{\sscp U}>2$ such that,
$$
f\in L^{p_{\sscp U}}\!\!\left(U\right)\mbox{ or either }
u\in W^{2,p_{\sscp U}}\!\!\left(U\right).
$$ 
Also, by setting $B_r(S)= \bigcup_{p\in S} B_r(p)$, we say that,
$$
f\in L^{p,{\sscp \rm loc}}\left(\om\setminus S\right)\mbox{ or either }
u\in W^{2, p, {\sscp \rm loc}}\left(\om\setminus S\right),\;\mbox{ for some }p>2,
$$ 
if for each $r>0$, there exists $p_r>2$ such that,
$$
f\in L^{p_{r}}\!\!\left(\om\setminus B_r(S)\right)\mbox{ or either }
u\in W^{2,p_{r}}\!\!\left(\om\setminus B_r(S)\right).
$$}
\edf

\bdf {\it
Let $f_{\pm}$ be two superharmonic functions in $\om$ taking the form \eqref{frep}, $u\in L^{1}_{\rm loc}(\om)$ and
$K e^{f+u}\in L^{1}_{\rm loc}(\om)$.
For any fixed and relatively compact Borel set  $E\Subset \om$, we define,
\beq\label{gammadefn}
\mathcal{K}_+(E;K_0)=\sup\limits_{U\subseteq E}\left \{\frac12\omega(U)+\int\limits_{U}[K-K_0]{\dsp e}^{f+u}dx\right \},
\eeq
where the supremum is taken over all Borel sets $U\subseteq E$.}
\edf

\bigskip

Finally we will need the following result about the local exponential integrability of $e^{f}$. Although similar 
exponential estimates for logarithmic potentials are well known, see \cite{Res3} or more recently \cite{bm} and \cite{Troy2}, 
it seems that the statement which is really needed here has been introduced only very recently in \cite{Ambrosio}.

\bpr\label{prexp}
Let $f_{\pm}$ be two superharmonic functions satisfying \eqref{frep} in $\om$ and
\beq\label{s2p}
S_{2\pi}=\{x\in \om\;:\;\omega_+(x)\geq 2\pi \}.
\eeq 
Then $S_{2\pi}$ is finite, $d_1=\frac{1}{4}\mbox{dist}(S_{2\pi},\pa\om)>0$ 
and we have:\\
$(i)$ $e^{-(f_--h_-)}\in L^{\ii}(\om)$ and $e^{(f_+-h_+)}\in 
L^{p_0,{\sscp \rm loc}}(\om\setminus S_{2\pi})$ for some $p_0>2$.\\
$(ii)$  If 
\beq\label{nocusp}
\fo\,x\in\om,\;\;\omega_+(x)<4\pi.
\eeq 
holds, then $e^{(f_+-h_+)}\in L^{q_0}(\om)$ for some $q_0> 1$.
\epr

Our main result is the following,

\bte\label{thm1} Let $K_0\geq 0$,  $f_{\pm}$ be two superharmonic functions taking the form \eqref{frep} and satisfying {\rm \rife{nocusp}}, $q_0>1$ and $p_0>2$ be defined as in Proposition \ref{prexp}
and $K\in  L^{n,{\sscp \rm loc}}_{\rm loc}(\om\setminus S_{2\pi}) \cap L_{\rm loc}^s\left(\om\right)$, 
for  some $n>\frac{2 p_0}{p_0-2}$ and some $s>\frac{q_0}{q_0-1}$.
Assume that, either,\\
$(i)$ $K e^{f+u}\in L^1\left(\om\right)$, where $u\in L^1\left(\om\right)$ is a solution of \eqref{1n} in the sense of distributions, 
or,\\
$(ii)$ $u\in W^{2,p,{\sscp \rm loc}}_{\rm loc}\left(\om\setminus S_{2\pi}\right)\cap W^{2,q}_{\rm loc}\left(\om\right)$, 
for some $p>2$ and some $q>1$, is a strong subsolution of \eqref{1n}, that is, 
\beq\label{1nn}
-\Delta u \leq   2K e^{f}\e{u}\;\; \mbox{ for a.a. }x\in\om.
\eeq

Then, if $(i)$ holds, we have
$u\in W^{2,p,{\sscp \rm loc}}_{\rm loc}\left(\om\setminus S_{2\pi}\right)\cap W^{2,q}_{\rm loc}\left(\om\right)\cap L^{\ii}_{\rm loc}(\om)$, 
for some $p>2$ and some $q\in(1,2)$, and in particular $u$ is a strong solution of \eqref{1n}, that is, it satisfies {\rm \rife{1nn}} with the equality sign.\\
Moreover,  in both cases, for any fixed simple and relatively compact subdomain $E\Subset \om$, we have $M(E)<+\ii$ and the following inequality holds: 
\beq\label{Alex.intron}
L^2(\pa E)\geq \left(4\pi-2\mathcal{K}_+(E;K_0)-K_0M(E)\right)M(E).
\eeq
The equality in \eqref{Alex.intron} is attained if and only if $u$ is a strong solution of \eqref{1n} in $E$ and,
\beq\label{sharp-last}
\e{f(z)+u(z)}=\frac{\tau^2\left| \Phi_0^{'}(z)( \Phi_0(z) )^{-\al}\right|^2}{(1+\frac{K_0\tau^2}{4(1-\al)^2} |\Phi_0(z)|^{2(1-\al)})^2},\;z\in E,
\eeq
for some $\tau\neq 0$,
where $\al=\frac{1}{4\pi}\omega_+(E)$ and $\Phi_0$ is a conformal map of $E$ onto the disk of unit radius, $|\Phi_0(z)|<1$ with
$\Phi_0(z_0)=0$, for some $z_0\in E$.\\
$(iii)$ If $\omega \perp e^{f+u} \mathcal{H}^2$, then the equality holds
if and only if, in addition to the
above conditions, one has $K\equiv K_0$ for $a.a.\,z\in E$ and $\omega=4\pi \al \dt_{z=z_0}$, that is, $f(z)=h(z)+4\pi\al G(z,z_0)=h(z)-2\al \log|\Phi_0(z)|$,
for some harmonic function $h$ in $E$.
\ete

\brm\label{rem1.5}{\it  Let us denote by $\mathcal{H}^2$ the 2-dimensional Hausdorff measure. By using the fact that $\omega_+\perp\omega_-$,
it is easy to check that if $\omega_-\perp{\dsp e}^{f+u}\mathcal{H}^2$, then we have,
$$
2\mathcal{K}_+(E;K_0)=\omega_+(E)+2\int\limits_{E}[K-K_0]^+{\dsp e}^{f+u}dx,
$$
where $K^+=\max\{K,0\}$, while in general the equality sign should be replaced by the inequality  sign.
}
\erm

As far as one is just concerned with the inequality and not with the characterization of the equality sign,
then, in case $(i)$ holds, \eqref{Alex.intron} holds under much weaker conditions. 
The proof of this fact is based on Theorem \ref{thm1} and on some results and arguments in \cite{bm} about the regularity properties of Liouville-type equations.
\bco\label{co17}
Let $K_0\geq 0$,  $f_{\pm}$ be two superharmonic functions taking the form \eqref{frep} and satisfying {\rm \rife{nocusp}}, 
and $K\in L^1(\om; e^{f+u}\mathcal{H}^2)$ where $u\in L^1\left(\om\right)$ is a solution of \eqref{1n} in the sense of distributions. 
Then:\\
$(i)$ $u\in W_{\rm loc}^{1,r}(\om)$ for any $r\in (1,2)$ and $e^{t|u|}\in L^{1}_{\rm loc}(\om)$ for any $t \geq 1$;\\
$(ii)$ \eqref{Alex.intron} holds.
\eco

Motivated by the study of a
cosmic string equation, in a recent paper \cite{BCast1} we derived \eqref{Alex.intron} in the easier situation where $\omega_-\equiv 0$ while $\omega_+$
is proportional to a Dirac delta.  The problem here is more subtle, and
the crux of the proof is to attach to each strong subsolution of \rife{1n} an auxiliary function (which we will denote by $\eta$)
which satisfies a Liouville type equation with
Dirichlet boundary condition on $E$, and which admits a suitable locally absolutely continuous weighted rearrangement (which we will denote by $\eta^*$).
The difficulty arises since, in view of the generality suggested by the geometric application, no
assumption is made about $\omega$, with the unique exception of the "no-cusp" hypothesis \rife{nocusp}. As a consequence,
the term $e^{f_+}$, which is part of the weight factor in the weighted rearrangement, can come with almost any kind of singularity.
In particular, the standard argument \cite{Tal}
yielding the absolute continuity of $\eta^*$, does not work in this case, neither in the slightly improved form used to handle conical singularities, see \cite{BCast1}.
We succeed in solving this problem by a careful decomposition of the
singular set of $\omega_+$, see the definition of $S_{2\pi}$ in \rife{s2p}.
The point is that $S_{2\pi}$ is finite in $\om$, while, locally in its complement, we come up with enough summability for $e^f$ to guarantee that
$\eta^*$ is absolutely continuous. This approach, recently pursued in \cite{Ambrosio} to prove a regularity result for a class of singular
surfaces introduced by Alexandrov \cite{Ale2}, motivates the peculiar notations introduced in Definition \ref{def12}.
In particular, the assumptions about $K$ and $u$, are essentially the
minimal requirements to match the regularity of $\eta$ as allowed by the properties of $f$ derived in this way. It is understood that 
the characterization of the equality case in this weak contest is new as well.\\
The proof of Theorem \ref{thm1} is split into four steps. In the first and second step we construct $\eta$, its weighted rearrangement $\eta^*$ and
prove that $\eta^*$ is locally absolutely continuous. Step four contains the discussion about the equality case.
Step three is the adaptation in our setting of the part of the
Bandle argument which is concerned with the derivation of a differential inequality and its consequences.\\

In the second part of this paper, and in the same spirit of \cite{Band0}, we will apply \eqref{Alex.intron} to derive a new proof of the Alexandrov isoperimetric
inequality for $K_0\geq 0$ on abstract surfaces of bounded curvature, see \rife{Alex.intro} in Theorem \ref{Alexandrovn}. We refer the reader 
to \cite{Band}, \cite{BaDo}, \cite{Bur}, \cite{Hub0}, \cite{Oss2} and the references
therein for a detailed exposition of the proof and of the interesting history of Alexandrov's inequality 
and to  \cite{topp}, \cite{topp1} for other more recent proofs. See also \cite{BLin3}, \cite{suz}. While in the above references one can find various proofs 
of the inequality \rife{Alex.intro}, we were not able to find a proof of the characterization of the equality case in the weak context pursued here, which 
seems therefore to be new even in the geometric setting.\\ 
Besides, to apply our estimates to this problem, 
we need to prove a seemingly new characterization of the structure of the metrics in local isothermal coordinates for certain classes of singular surfaces,
see Theorem \ref{thm2}. This intermediate result can also be seen as another result in the description
of the regularity properties of isothermal coordinates systems on Alexandrov's surfaces of bounded curvature recently pursued in \cite{Ambrosio}. 
Finally, some explicit examples are discussed to illustrate these results, including the isoperimetric inequality \rife{Alex.intro} 
on various singular surfaces homeomorphic to the 2-sphere.

\bigskip

We conclude this introduction with a remark about the case where $E$ is not simple but just regular, that is, the possibility that $E$ could be connected 
but not simply connected.
\brm\label{multi}{\it
If $\om$ is simply connected and the assumptions of Theorem \ref{thm1} are satisfied, 
but if the set $E\Subset \om$ is just assumed to be regular, then it is straightforward to check that
our proof yields the following inequality,
\beq\label{alelast}
L^2(\pa E)> \left(4\pi-2\mathcal{K}_+(E_s;K_0)-K_0M(E)\right)M(E),
\eeq
where $E_s$ is the interior of the closure of the union of $E$ with the bounded components of $\RR^2\setminus E$ (the "holes" of $E$) , 
which we denote by $(E)_B$, that is, 
$$
E_s=\accentset{\circ}{\ov{E \cup (E)_B}}.
$$
In other words we still have an inequality of the form 
\eqref{Alex.intron}, but we have a worse isoperimetric ratio, which is essentially obtained by subtracting the terms of the total curvature relative to the 
"holes" of $E$. This is not a technical point, and in fact it is possible to construct counterexamples to the inequality 
where these terms are omitted, see for example p.14 in \cite{BLin3}. The proof of this inequality is really the same as that of Theorem \ref{thm1}, 
but for the fact that in \eqref{Huber-subn} and in \eqref{sharp4} below we use the Huber inequality \eqref{Hubineq-sing-1} for the non contractible 
domain $E$. In particular this is also why we obtain the strict inequality in this case. It is straightforward to check that if the 
assumptions of Corollary \ref{co17} are satisfied, then \eqref{alelast} holds with the $\geq$ sign replacing the strict inequality.}
\erm

This paper is organized as follows. In section \ref{sec2} we prove Proposition \ref{prexp} and discuss the Huber's inequality. In section \ref{sec3} we prove 
Theorem \ref{thm1} and Corollary \ref{co17}. Sections \ref{sec4} and \ref{sec5} are devoted to the discussion of the Alexandrov's isoperimetric inequality and related examples. 

\bigskip
\bigskip

\section{Preliminary estimates: exponential summability of subharmonic functions and Huber's inequality.}\label{sec2}

The local exponential integrability of $e^{f_+}$ as claimed in Proposition \ref{prexp} is not new, see \cite{Ambrosio}.
We provide the proof of Proposition \ref{prexp} for the sake of completeness.\\

{\it The Proof of Proposition \ref{prexp}.}\\
We will denote by $d_{\sscp \om}$ the diameter of $\om$. Clearly $S_{2 \pi}$ is finite since $\omega_+$ is finite, whence obviously 
$\mbox{dist}(S_{2\pi},\pa\om)>0$. Let $d_1=\frac14\mbox{dist}(S_{2\pi},\pa\om)$ and
let us set $\om_d=\{x\in \om: \mbox{dist}(x,\pa\om)<d\}$. Then
$\omega_{+}(\om_{d})\searrow 0^+$, as $d \searrow 0^+$, whence there exists $d_0>0$ such that
$\omega^{+}(\om_{d})<\frac{\pi}{2}$, for each $d<4d_0$. We choose $d_0$ possibly smaller to satisfy $4d_0<d_1$.
It is not difficult to see that there exists $C_0>0$ such that,
$$
(f_+(x)-h_+(x))-C_0\leq w_0 (x) := \frac{1}{2\pi} \int_{\om_{2d_0}} \log \left( \frac{d_\om}{|x-y|} \right) \, d \omega_{+},\;\fo\;x\in \om_{d_0}.
$$

By the Jensen's inequality and the Fubini-Tonelli Theorem we can estimate,
$$
\int_{\om_{d_0}} \exp \left(\frac{ 3\pi  w_0}{\omega_{+} (\om_{2 d_0})} \right) d x  \leq $$
$$
\int_{\om_{d_0}} d x \int_{\om_{2 d_0}} \left( \frac{d_{\sscp \om}}{|x-y|} \right)^{\frac{3}{2}} \frac{d \omega_{+} (y)}{\omega_{+} (B_{2 d_0})} =
\int_{\om_{2 d_0}} \frac{d \omega_{+} (y)}{\omega_{+} (\om_{2 d_0})} \int_{\om_{d_0}}
\left( \frac{d_{\sscp \om}}{|x-y|} \right)^{\frac{3}{2}} d x\leq
$$
$$\int_{\om_{2 d_0}} \frac{d \omega_{+} (y)}{\omega_{+} (\om_{2 d_0})} \int_{B_{d_{\sscp \om}} (y)}
\left( \frac{d_{\sscp \om}}{|x-y|} \right)^{\frac{3}{2}} d x = \pi (2d_{\sscp \om})^2,
$$
where we used the fact that $\om_{d_0} \subset B_{d_{\sscp \om}} (y)$. This inequality shows that $e^{(f_+-h_+)}\in L^{6}(\om_{d_0})$.\\

$(i)$ Since $-(f_--h_-)$ is negative, then $e^{-(f_--h_-)}\in L^{\ii}(\om)$.
Let $\om_0 ={ \{\om\setminus \om_{\frac{d_0}{2}}\}\setminus B_{r}(S_{2 \pi})}$, with $0<r<d_1$,
and let us fix $x_0 \in \ov{\om_0}$. Since $\omega_+(x_0)<2\pi$,
then we can find $\eps>0$ such that there exists $R>0$ depending on $x_0$ and $\eps$, such that the ball centred at $x_0$, $B_{2R}:= B_{2R}(x_0)$,
satisfies $B_{2R} \Subset\{\om\setminus \om_{d_0}\} \setminus S_{2 \pi}$
and $\omega_+(B_{2R}) \leq 2\pi - 2\eps$.
As above there exists $C>0$ such that,
$$
(f_+(x)-h_+(x))-C\leq w (x) := \frac{1}{2\pi} \int_{B_{2R}} \log \left( \frac{4R}{|x-y|} \right) \, d \omega_{+},\;\fo\;x\in D_{x_0}\equiv D_0:=B_{R}(x_0),
$$
and for any $\delta < 4\pi$ we can estimate,

$$
\int_{D_0} \exp \left(\frac{(4\pi - \delta) w}{\omega_{+} (B_{2R})} \right) d x \leq $$
$$\int_{D_0} d x \int_{B_{2R}} \left( \frac{d_{\sscp \om}}{|x-y|} \right)^{2 - \frac{\delta}{2\pi}} \frac{d \omega_{+} (y)}{\omega_{+} (B_{2R})} =
\int_{B_{2R}} \frac{d \omega_{+} (y)}{\omega_{+} (B_{2R})} \int_{D_0} \left( \frac{d_{\sscp \om}}{|x-y|} \right)^{2 - \frac{\delta}{2\pi}} d x\leq $$
$$\int_{B_{2R}} \frac{d \omega_{+} (y)}{\omega_{+} (B_{2R})} \int_{B_{d_{\sscp \om}} (y)} \left( \frac{d_{\sscp \om}}{|x-y|} \right)^{2 - \frac{\delta}{2\pi}} d x =
\frac{(2\pi d_{\sscp \om})^2}{\delta}.
$$
Therefore, in particular by choosing $\delta < \eps$,
we see that $p_{\sscp D_0} := \frac{(4\pi - \delta)}{\omega_{+} (D_0)}> \frac{(4\pi - \eps)}{\omega_{+} (B_{2R})}> 2$ so that
$e^{(f_+-h_+)} \in L^{p_{\sscp D_0}} (D_0),$ for some $p_{\sscp D_0} >2$ depending on $x_0$ and $R$.\\
At this point we define
$\mathcal{B}=\bigcup_{x\in \ov{\om_0}}D_x$, where each $D_x$, constructed as above, comes with its own $p_{\sscp D_{x}}>2$.
Clearly $\mathcal{B}$ is an open cover of $\ov{\om_0}$, and since $\ov{\om_0}$ is compact, then we can extract a finite cover $D_{x_j}$, $j=1,\dots,N$,
and set $p_{\sscp U} := \min\left\{6, \min\limits_{j=1,..,N}p_{\sscp D_{x_j}}\right\}$.
Therefore $e^{(f_+-h_+)} \in L^{p_{\sscp U}} (\om \setminus B_r(S_{2\pi}))$, for some $p_{\sscp U}>2$, which proves $(i)$.\\

$(ii)$ Let us define $\om_1={\om\setminus \om_{\frac{d_0}{2}}}$. We use \eqref{nocusp},
as in the proof of $(i)$ to conclude that $e^{(f_+-h_+)}\in L^{q}(\om_1)$ for some
$q > 1$. Therefore we find $e^{(f_+-h_+)}\in L^{q_0}(\om)$ where $q_0=\min\{6,q\}>1$, as claimed.
\fineproof

\bigskip
\bigskip

Next we present the well known Huber's inequality \cite{Hub} as well as a generalization of it
suitable to be applied to regular (whence in particular non simply connected) domains.

\bte[The Huber inequality, \cite{Hub}] \label{hub-teo}Let $\om\subset \RR^2$ be open and bounded and
$E\Subset \om$ be a simple and relatively compact subset. Let $f$ be the difference of two superharmonic functions in $\om$
taking the form \eqref{frep}. Then it holds,
\beq\label{Hubineq-sing}
\left(\,\int\limits_{\pa E}{\dsp e}^{\frac{f}{2}}d\ell\right)^2
\geq \left(4\pi -\omega_+(E)\right)\int\limits_{E}{\dsp e}^{f}dx.
\eeq
The equality holds in \eqref{Hubineq-sing} if and only if, in complex notations, $f(z)=c+2\log\left|\Phi^{'}(z)( \Phi(z)^{-\dsp \al_{\sscp E}} )\right|$ where
$\al_{\sscp E}=\frac{1}{4\pi}\omega_+(E)$ and $\Phi$ is a conformal map of $E$ onto the disk of unitary radius $|w|=|\Phi(z)|<1$ with $\Phi(z_0)=0$ for
some $z_0\in E$.
\ete

\bigskip

We will need the following generalization of the Huber's result.
\bte\label{hub-teo-2} Let $\om\subset \RR^2$ be open and bounded and
$E\Subset \om$ be a simple and relatively compact subset. Let $f$ be the difference of two superharmonic functions in $\om$
taking the form \eqref{frep}. If $U\subseteq E$ is a regular domain, then it holds,

\beq\label{Hubineq-sing-1}
\left(\;\int\limits_{\pa U}{\dsp e}^{\frac{f}{2}}d\ell\right)^2
\geq \left(4\pi -\omega_+(E)\right)\int\limits_{U}{\dsp e}^{f}dx.
\eeq
In particular, if $U$ is not simply connected, then the inequality is strict.
\ete
\proof
In view of Theorem \ref{hub-teo} we are left to discuss the cases where $U$ is not simply connected
and prove in particular that in all those cases the inequality is strict. Obviously the inequality is trivially satisfied if $\omega_+(E)\geq 4\pi$, whence
we assume w.l.o.g. that $\omega_+(E)< 4\pi$.
Let us assume for the moment that $U=U_1\setminus\ov{U_0}$ for a pair of
simple domains such that $U_0\Subset U_1$ and
$\pa U = \pa U_1 \cup \pa U_0$. So $U_1=U\cup\ov{U_0}$ and in this case, by assumption we have $E=U_1$ and in particular $\omega_+(U_0)< \omega_+(U_1)<4\pi$.
For any domain $U\subset\rdue$, let us set
$$
\ell(\pa U)=\;\int\limits_{\pa U} {\dsp e}^{\frac{f}{2}}d\ell,\qquad
M(U)=\int\limits_{U} {\dsp e}^{f}dx.
$$
Thus we may use \rife{Hubineq-sing} to obtain
$$
\ell^2(\pa U)=\ell^2(\pa U_1 \cup \pa U_0)>
\ell^2(\pa U_1)+\ell^2(\pa U_0)\geq
$$
$$
\left(4\pi - \omega_+(U_1) \right)M(U_1)+
\left(4\pi -\omega_+(U_0)\right)M(U_0)>
$$
$$
\left(4\pi-\omega_+(U_1) \right) M(U_1)>\left(4\pi -\omega_+(U_1)\right) M(U),
$$
which is \eqref{Hubineq-sing-1} in this particular case.
The case where $\RR^2\setminus U$ has finitely many bounded components readily follows by an
induction argument on the number of "holes" of $U$. Obviously the inequality is always strict whenever $U$ is not simply connected.\fineproof

\bigskip

\section{The proof of Theorem \ref{thm1} and Corollary \ref{co17}.}\label{sec3}
This section is devoted to the proof of Theorem \ref{thm1} and Corollary \ref{co17}.

\bigskip

{\bf The proof of Theorem \ref{thm1}.}\\
Once the result has been established for $K_0\neq 0$, then the case $K_0=0$ is worked out by
an elementary limiting argument, which is why we will just discuss the case $K_0>0$.\\
We recall that by assumption $\frac{n p_0}{n+p_0}>2$ and $\frac{s q_0}{s+q_0}>1$.
First of all, we have the following,
\ble\label{lem0308}
$(a)$ If $(i)$ holds and if $K\in L_{\rm loc}^s\left(\om\right)$ for some $s>\frac{q_0}{q_0-1}$, then,
$$
Ke^{f+u}\in L^{r}_{\rm loc}\left(\om\right),
$$
and
$u\in L^{\ii}_{\rm loc}(\om)\cap W^{2,r}_{\rm loc}\left(\om\right)$, 
for any $1\leq  r\leq \frac{s q_0}{s+q_0}$. In particular $u$ is a strong solution of \eqref{1n}.\\

$(b)$  If $(i)$ holds and if $K\in  L^{n,{\sscp \rm loc}}_{\rm loc}(\om\setminus S_{2\pi}) \cap L_{\rm loc}^s\left(\om\right)$ 
for  some $n>\frac{2 p_0}{p_0-2}$ and some $s>\frac{q_0}{q_0-1}$, then,
\beq\label{300616.1}
Ke^{f+u}\in L^{k,{\sscp \rm loc}}_{\rm loc}\left(\om\setminus S_{2\pi}\right)\cap L^{r}_{\rm loc}\left(\om\right),
\eeq
and
$u\in W^{2,k,{\sscp \rm loc}}_{\rm loc}\left(\om\setminus S_{2\pi}\right)\cap W^{2,r}_{\rm loc}\left(\om\right)\cap L^{\ii}_{\rm loc}(\om)$, 
for any $2< k\leq \frac{n p_0}{n+p_0}$ and $1\leq  r\leq \frac{s q_0}{s+q_0}$. In particular $u$ is a strong solution of \eqref{1n}.
\ele
\proof $(a)$
By assumption we have $\frac{sq_0}{s+q_0}>1$ and  then, in view of Proposition \ref{prexp}, 
we also have $Ke^{f}\in L_{\rm loc}^{q}(\om)$, $\;\fo\,1< q\leq \frac{sq_0}{s+q_0}$. 
On the other hand, since $Ke^{f+u}\in L^1(\om)$, then, by Remark 2 in \cite{bm}, we have $e^{|u|}\in L^{k}_{\rm loc}(\om)$ for any $k>0$, and therefore 
in particular $e^u\in L^{q^{'}}(\om)$, where $q^{'}=\frac{q}{q-1}<+\ii$. Thus we can apply another result in \cite{bm} (see Remark 5 in \cite{bm}), 
which yields $u\in L^{\ii}_{\rm loc}(\om)$. So, 
by standard elliptic estimates, we conclude also that $u\in W_{\rm loc}^{2,r}(\om)$ and
in particular that $u$ is a strong solution of \eqref{1n}.\\ 
$(b)$  
Next, let us fix a compact set $U \subset \om\setminus S_{2\pi}$ and observe that, by assumption, $K\in L^{n}(U)$ 
for some $n>\frac{2 p_{\sscp U}}{p_{\sscp U}-2}$. Therefore $\frac{n  p_{\sscp U}}{n+ p_{\sscp U}}>2$ and then, in view of Proposition \ref{prexp}, 
we also have $Ke^{f}\in L^{p}(U),\;\fo\,2< p\leq \frac{n p_{\sscp U}}{n+p_{\sscp U}}$. Since $U$ is arbitrary, then we conclude that 
$Ke^{f+u}\in L_{\rm loc}^{k,{\sscp \rm loc}}(\om\setminus S_{2\pi})\cap L_{\rm loc}^{r}(\om)$, for any 
$2< k\leq \frac{n p_{\sscp U}}{n+p_{\sscp U}}$ and $1<  r\leq \frac{s q_0}{s+q_0}$. As above,
by standard elliptic estimates, we conclude also that $u\in W_{\rm loc}^{2,k,{\sscp \rm loc}}(\om\setminus S_{2\pi})\cap W_{\rm loc}^{2,r}(\om)$ and
in particular that $u$ is a strong solution of \eqref{1n}.\finedim

\bigskip
\bigskip

Lemma \ref{lem0308} shows that if $(i)$ holds, 
then $u\in W^{2,p,{\sscp \rm loc}}_{\rm loc}\left(\om\setminus S_{2\pi}\right)\cap W^{2,q}_{\rm loc}\left(\om\right)\cap L^{\ii}_{\rm loc}(\om)$,  for some 
$p>2$ and $q>1$, and moreover that $u$ is a strong solution of \rife{1n}. 
Whence we are reduced to the analysis of the case where 
$u\in W^{2,p,{\sscp \rm loc}}_{\rm loc}\left(\om\setminus S_{2\pi}\right)\cap W^{2,q}_{\rm loc}\left(\om\right)\cap L^{\ii}_{\rm loc}(\om)$,  for some 
$p>2$ and $q>1$, satisfies \eqref{1nn}. In particular, in the rest of the proof, we will use the fact that, by the Sobolev embedding Theorem, $u\in C^{0}_{\rm loc}(\om)$. 
Clearly, in view of  \rife{300616.1}, $M(E)$ is finite. We divide the proof into four steps.\\
{\bf Step 1.}\\
Since $E\Subset\om$ is relatively compact and simple,  then we can find
an open, simply connected,  relatively compact and smooth domain $\om_0$ such that,
$$
E \Subset \om_0\Subset \om.
$$
Since $S_{2\pi}$ is finite and since $\omega_{\pm}(\om_0)<+\ii$, then we can choose $\om_0$ such that, for some $N\in\NN$,
\beq\label{spi}
S^{\sscp 0}_{2\pi}:=S_{2\pi}\cap \om_0=\{q_1,\dots,q_{\sscp N}\}\subset \om_0\mbox{ and }\pa\om_0\cap S_{2\pi}=\emptyset.
\eeq

Clearly, in view of \rife{1nn}, we have,
\beq\label{13n}
-\Delta u \leq 2 K e^{f} \e{u} = 2[K-K_0] e^{f}\e{u}+2K_0e^{f} \e{u}\;\;\mbox{ for a.a. }x\in \om_0.
\eeq
Next, let us define,
\beq\label{sharp3}
\phi(x):=-\Delta u - 2[K-K_0] e^{f}\e{u}-2K_0 e^{f}\e{u},\quad  x\in \om_0.
\eeq

Since $u\in W^{2,p,{\sscp \rm loc}}_{\rm loc}\left(\om\setminus S_{2\pi}\right)\cap W^{2,q}\left(\om\right)\cap C^{0}_{\rm loc}(\om)$,
for some $p>2$ and some $q> 1$, and in view of \rife{300616.1} and of Proposition \ref{prexp}, we see from \rife{13n} that,
$$
\phi(x)\leq 0, \;\;\;\mbox{ for a.a. }x\in  \om_0\;\;\mbox{ and }\;\;
\phi\in L^{p, {\rm \sscp loc}}_{\rm loc}(\om_0\setminus S^{\sscp 0}_{2\pi})\cap L^{q}(\om_0),
$$
for some $p>2$ and some $q> 1$. Therefore, in view of Theorem 9.15, Corollary 9.18 and Lemma 9.17 in
\cite{Gilb} we see that the linear problem,
\beq\label{0403.10nn}
\Delta w = \phi\;\; \mbox{in}\;\;\om_0,
\qquad w = 0 \;\; \mbox{on}\;\; \pa \om_0,
\eeq
admits a unique strong solution $w \in W^{2,p, {\rm \sscp loc}}_{\rm loc}(\om_0\setminus S^{\sscp 0}_{2\pi})\cap W^{2,q}(\om_0)\cap C^{0}(\,\ov{\om_0}\,)$,
for some $p>2$ and some $q>1$. Obviously $w$ is superharmonic (see \cite{Gilb} \S 2.8 and Ex. 2.7, 2.8).\\

Next let $f_1$ be the Perron's (see \S 2.8 in \cite{Gilb}) solution of
$\Delta f_1=0$ in $ E$, $f_1=-u$ on $\pa E$. Since $u\in C^{0}(\overline{ E})$, then $f_1$ is well
defined and continuous up to the boundary (see \S 2.8 in \cite{Gilb}). Let us also define
$f_{2}$ to be the unique $W^{2,p, {\rm \sscp loc}}_{\rm loc}(\om_0\setminus S^{\sscp 0}_{2\pi})\cap W^{2,q}(\om_0)\cap C^{0}(\,\ov{\om_0}\,)$
(for some $p>2$ and some $q>1$) solution of the linear problem,
$$
-\Delta f_{2} = 2[K-K_0] e^{f}\e{u} \;\; \mbox{in} \;\; \om_0,
\qquad f_{2}= 0 \;\; \mbox{on}\;\; \pa \om_0.
$$

With these definitions, we may finally set $\eta=u+w+f_1-f_{2}$.
Then, we see that $\eta\in  W^{2,p, {\rm \sscp loc}}_{\rm loc}(E\setminus S^{\sscp 0}_{2\pi})\cap W^{2,q}(E)\cap C^{0}(\,\ov{E}\,)$
for some $p>2$ and some $q>1$ and satisfies,
\beq\label{eqeta-subn}
-\Delta \eta= 2K_0 \e{\psi}\e{\eta}\;\; \mbox{ for a.a. }x\in   E,
\quad \eta=0\;\; \mbox{on}\;\; \pa E,
\eeq
where
\beq\label{subh-subn}
\psi=f_{+}+f_{2}-f_{-}-w-f_1.
\eeq

\bigskip

By the Sobolev embedding Theorem we conclude that,
\beq\label{reg.1b-subn}
\eta_+\in C^{1}_{\rm loc}(E\setminus S^{\sscp 0}_{2\pi}).
\eeq

Since $\eta \in W^{2,q}(E)$, for some $q >1$, then by using the Sobolev embedding once more we see that $\eta \in W^{1,2}(E) \cap C^{0} (E)$.
Then by the maximum principle for weak solutions (see for example Theorem 8.1 in \cite{Gilb}) we deduce that $\eta\geq 0$. In particular,
by the strong maximum principle for weak supersolutions (see for example Theorem 8.18 in \cite{Gilb}) we also check that $\eta$ is strictly positive in $E$.
In particular, we conclude that,
\beq\label{2.71bisn}
\eta(x)>0\quad \fo\, x\in E\quad\mbox{and}\quad \eta(x)=0\;\iff\;x\in\pa E.
\eeq

\bigskip

{\bf Step 2.}\\
Let us set $t_{+}=\max\limits_{\ov{ E}} \eta$,
$$
d\tau={\dsp e}^{\psi}dx,\qquad \sigma={\dsp e}^{\frac{\psi}{2}}d\ell,
$$
and let us define,
$$
\om(t)=\{x\in E\,|\,\eta (x)>t\},\;\; t\in  [0,t_+),\quad
\Gamma(t)=\{x\in E\,|\,\eta(x)=t\},\;\; t\in  [0,t_+],
$$
and
$$
\mu(t)=\int\limits_{\om(t)}d\tau.
$$

Since $\eta$ satisfy \rife{eqeta-subn}, then $\Gamma(t)$ has null measure,  whence  we
conclude that $\mu$ is continuous. Moreover, in view of \rife{2.71bisn}, we find that,
\beq\label{reg.2-subn}
\om(0)= E,\qquad \Gamma(0)=\pa E,\qquad \mu(0)=\int\limits_{ E}d\tau.
\eeq

Clearly we can extend $\mu$ on $[0,t_+]$
by setting $\mu(t_+)=\lim\limits_{t\nearrow t_+}\mu(t)=0^{+}$,
whence $\mu\in C^{0}([0,t_+])$. Next, by using \rife{eqeta-subn} once more, it is not difficult to see that
the 2-dimensional measure of the set $\{x\in E\::\: \nabla \eta(x)=0\}$ vanishes.  Therefore, by a well known consequence of the co-area
formula (see for example \cite{bz} p.158) and of the Sard's Lemma for Sobolev functions \cite{dep} (here we use also \eqref{spi}), we see that
\beq\label{070214.1-subn}
\frac{d\mu(t)}{dt}=-\int\limits_{\Gamma(t)}\frac{{\dsp e}^{\psi}}{|\nabla \eta|}\,d\ell,
\eeq
$\mbox{for a.a. }t\in [0,t_+]$.\\
At this point, for any $s\in[0,\mu(0))\equiv [\mu(t_+),\mu(0))$, we introduce a weighted
rearrangement of $\eta$,
\beq\label{etadef-subn}
\eta^*(s)= |\{t\in [0,t_+]:\mu(t)>s\}|,
\eeq
where $|U|$ denotes the Lebesgue measure of a Borel set $U\subset \R$. By setting
$\eta^*(\mu(0))=0$, then $\eta^*\in C^{0}([0,\mu(0)])$ is
the inverse of $\mu$ on $[0,t_+]$ and
coincides with the distribution function of $\mu$. Actually $\eta^*$ is strictly decreasing, whence
differentiable almost everywhere. A crucial point at this stage is to prove
that $\eta^*$ is not just continuous and differentiable almost everywhere, but also locally absolutely continuous.
It turns out that in fact it is locally Lipschitz in $(0,\mu(0))$ as shown in the following Lemma:

\ble\label{lem-liploc}
For any  $0<\overline{a}\leq a<b\leq \overline{b}<\mu(0)$, there exist
$\ov{C}=\ov{C}(\overline{a},\overline{b},S_{2\pi},\mathcal{K}_+(E;K_0))>0$ such that,
\beq\label{liploc}
\eta^*(a)-\eta^*(b)\leq \ov{C} (b-a).
\eeq
\ele
\proof
In view of \eqref{spi} and \eqref{reg.1b-subn}, we see that $|\nabla \eta| \leq C_{U}$ on any $U \Subset E \setminus S^0_{2\pi}$.
Let us then set $t_i = \eta (x_i)$,
$x_i \in S_{2\pi}$, $i=1,..,m$, with $m\leq N$, and $t_0 = \eta^* (a)$ and $t_{m+1} = \eta^* (b)$.
For any
$$
\eps< \min\left\{\frac{|\eta^*(a)-\eta^*(b)|}{4(m+1)},\; \frac14 \min\limits_{i=0,\dots,m}\{t_{i+1}-t_i\}\right\},
$$
we can find $\delta = \delta_\eps$ such that
$\eta^{-1} [t_i + \eps , t_{i+1}- \eps ] \cap B_{\delta} (S_{2\pi}) = \emptyset$ for any $i=0,..,m$,
where $B_{\delta} (S_{2\pi})$ is a $\delta$-neighbourhood of the set $S_{2\pi}$. Therefore, in particular, we can find $C_\eps>0$ such that
$|\nabla \eta(x)| \leq C_{\eps}$,  $\fo\;x\in\eta^{-1} [t_i + \eps , t_{i+1}- \eps]$.
At this point, since $K_0>0$, then  we can assume w.l.o.g. that $2\gamma_{ E}(K_0):=4\pi-2\mathcal{K}_+(E;K_0)>0$
(otherwise $4\pi-2\mathcal{K}_+(E;K_0)-K_0 M(E)<0$
and \eqref{Alex.intron} would be trivially satisfied). Therefore we can use the
coarea formula (see \cite{bz} p. 158) and the Huber's isoperimetric inequality \eqref{Hubineq-sing-1}, to conclude that,

$$
b-a=\mu(\eta^*(b))-\mu(\eta^*(a))=
\int\limits_{\eta>\eta^*(b)} d\tau-\int\limits_{\eta>\eta^*(a)} d\tau=
\int\limits_{\eta^*(b)<\eta\leq\eta^*(a)} \hspace{-.5cm}d\tau\geq
$$
$$
\int\limits_{\eta^*(b)<\eta<\eta^*(a)} \hspace{-.5cm}d\tau=
\int\limits_{\eta^*(b)}^{\eta^*(a)}
\left(\, \int\limits_{\Gamma(t)}\frac{d\sigma}{|\nabla \eta|}\right)dt=
\sum_{i=0}^{m} \int\limits^{t_{i+1}}_{t_i}\left(\, \int\limits_{\Gamma(t)}\frac{d\sigma}{|\nabla \eta|}\right)dt\geq
$$
$$
\sum_{i=0}^{m} \int\limits^{t_{i+1} - \eps}_{t_i + \eps}
\left(\, \int\limits_{\Gamma(t)}\frac{d\sigma}{|\nabla \eta|}\right)dt \geq
\frac{1}{C_\eps}  \sum_{i=0}^{m} \int\limits^{t_{i+1} - \eps}_{t_i + \eps}
\left(\, \int\limits_{\Gamma(t)}d\sigma\right)dt\geq
$$
$$
\frac{\sqrt{2\gamma_{ E}(K_0)}}{C_\eps}\sum_{i=0}^{m} \int\limits^{t_{i+1} - \eps}_{t_i + \eps}  \sqrt{\left(\, \int\limits_{\om(t)}d\tau\right)}\geq
\frac{\sqrt{2\gamma_{ E}(K_0)}}{C_\eps}\sqrt{\left(\, \int\limits_{\om(\eta^*(\ov{b}))}d\tau\right)}\sum_{i=0}^{m} \int\limits^{t_{i+1} - \eps}_{t_i + \eps} dt=
$$
$$
C(\ov{a},\ov{b},S_{2\pi},\mathcal{K}_+(E;K_0))|\eta^*(a)-\eta^*(b) - 2 (m+1) \eps|\geq\frac14 C|\eta^*(a)-\eta^*(b)|,
$$
for a strictly positive constant $C$ depending on $\ov{a},\ov{b},S_{2\pi},\mathcal{K}_+(E;K_0)$, as claimed.
\fineproof

\bigskip

{\bf Step 3.}\\
In view of \rife{070214.1-subn} we obtain,
\beq\label{070214.1-sub-1n}
\frac{d\eta^*(s)}{ds}=-
\left(\,\int\limits_{\Gamma(\eta^*(s))}\frac{{\dsp e}^{\psi}}{|\nabla \eta|}\,d\ell\right)^{-1},
\eeq
for any $s\in I^{*}$, where $[0,\mu(0)]\setminus I^{*}$ is a set of null measure and, by setting
$I:=\eta^{*}(I^{*})$, then $\mu(I)=I^*$. Next, let us define,
$$
F(s)=2K_0\int\limits_{\om(\eta^{*}(s))}{\dsp e}^{\eta}d\tau,\quad s\in[0,\mu(0)],
$$
where,
\beq\label{fin-sub1n}
F(\mu(0))=2K_0\int\limits_{ E}{\dsp e}^{\eta}d\tau=2K_0 M( E),
\eeq
and we have set,
\beq\label{fin-sub0n}
F(0)=\lim\limits_{s\searrow 0^+}F(s)=0^{+}.
\eeq
Clearly $F(s)$ is strictly increasing and
continuous on  $[0,\mu(0)]$ and in particular locally Lipschitz in $(0,\mu(0))$,
since in fact it satisfies,
$$
|F(s)-F(s_0)|\leq C |\mu(\eta^{*}(s))-\mu(\eta^{*}(s_0))|= C|s-s_0|,\quad \fo\,
0=\mu(t_+)<s_0<s<\mu(0),
$$
for a suitable constant $C>0$. In particular it holds,
$$
\int\limits_{\om(\eta^*(s))}{\dsp e}^{\dsp u}d\tau=
\int\limits_{0}^{s}{\dsp e}^{\dsp\eta^*(\lm)}d\lm,\;\;\;\fo\,s\in[0,\mu(0)],
$$
so that,
\beq\label{F2n}
\frac{dF(s)}{ds}= 2K_0 e^{\eta^*(s)},\quad
\frac{d^2 F (s)}{ds^2}=2K_0\frac{d\eta^*(s)}{ds}\,{\dsp e}^{\dsp \eta^*(s)}=
\frac{d\eta^*(s)}{ds}\frac{d F(s)}{ds},\quad\fo\,s\in I^*.
\eeq
We remark that since $\eta^*(s)$ is differentiable almost everywhere,
then the formula for the first derivative of $F(s)$ shows that in fact $\frac{dF(s)}{ds}$ is differentiable almost everywhere as well.\\
For any $s\in I^*$ the Cauchy-Schwartz inequality yields,
\beq\label{281109.1.1n}
\left(\,\int\limits_{\Gamma(\eta^*(s))} d\sigma\right)^2\leq
\left(\,\int\limits_{\Gamma(\eta^*(s))}\frac{{\dsp e}^{\psi}}{|\nabla \eta|}\,d\ell\right)
\left(\,\int\limits_{\Gamma(\eta^*(s))}|\nabla \eta|d\ell\right)=
\eeq
$$
\left(-\frac{d\eta^*(s)}{ds}\right)^{-1}
\left(\,\int\limits_{\Gamma(\eta^*(s))}\left(-\frac{\pa \eta}{\pa \nu_+}\right)d\ell\right),
$$
where $\nu_+=\frac{\nabla \eta}{|\nabla \eta|}$ is the exterior unit normal to $\om(\eta^*(s))$
and we have used
\rife{070214.1-sub-1n}. Obviously, we can assume w.l.o.g. that $\eta^{-1}(S_{2\pi}\cap E)\notin I$, so that,
since $\eta$ satisfies \rife{reg.1b-subn}, then \eqref{eqeta-subn} readily implies that,
$$
\int\limits_{\Gamma(\eta^*(s))}\left(-\frac{\pa \eta}{\pa \nu_+}\right)d\ell=
\int\limits_{\om(\eta^*(s))}2K_0 \e{\eta}d\tau,
$$
for any $s\in I^*$. Therefore, in particular we deduce that,
$$
\int\limits_{\Gamma(\eta^*(s)))}\left(-\frac{\pa \eta}{\pa \nu_+}\right)d\ell=
\int\limits_{\om(\eta^*(s))}2K_0 \e{\eta}d\tau=F(s),
$$

for any $s\in I^*$.
Plugging this identity in \rife{281109.1.1n} we find,
\beq\label{iso-1-subn}
\left(\,\int\limits_{\Gamma(\eta^*(s))} d\sigma\right)^2
\leq \left(-\frac{d\eta^*(s)}{ds}\right)^{-1}  F(s),
\eeq
for any $s\in I^*$.
Clearly, in view of \eqref{subh-subn}, we have,
\beq\label{sharp2}
-\Delta \psi=\omega_{+}-\omega_{-}+\phi-\Delta f_{2}\leq \omega+2[K-K_0] e^{f}\e{u},
\eeq
whence
\beq\label{sharp1}
\sup\limits_{U\subset E}\left \{\int\limits_{U}(-\Delta \psi)\right \}  \leq 2\mathcal{K}_+(E;K_0),
\eeq
and we can apply the generalized Huber's inequality \eqref{Hubineq-sing-1}, to conclude that,
\beq\label{Huber-subn}
\left(\,\int\limits_{\Gamma(\eta^*(s))} d\sigma\right)^2\geq
[4\pi-2\mathcal{K}_+(E;K_0)]\mu(\eta^*(s))\equiv
[4\pi-2\mathcal{K}_+(E;K_0)]s,
\eeq
$\mbox{for any}\;s\in I^*\cap (0,\mu(0))$.

\brm\label{rmmass}{\it
If $4\pi -2\mathcal{K}(E;K_0^+)<0$, then \eqref{Huber-subn} trivially satisfied.}
\erm

To simplify the exposition let us set,
$$
2\gamma_{ E}(K_0)=4\pi-2\mathcal{K}_+(E;K_0).
$$

Hence, substituting \rife{Huber-subn} in \rife{iso-1-subn}, we obtain,
$$
2\gamma_{ E}(K_0)s\leq\left(-\frac{d\eta^*(s)}{ds}\right)^{-1}
F(s)
,\,\mbox{for any}\;s\in I^*\cap (0,\mu(0)).
$$
So, multiplying by $\frac{dF(s)}{ds}\left(-\frac{d\eta^*(s)}{ds}\right)$, we come up with the inequality,
$$
2\frac{dF(s)}{ds}\left(\frac{d\eta^*(s)}{ds}\right)\gamma_{ E}(K_0)s
+\frac{dF(s)}{ds}F(s)\geq 0,\,
\mbox{for any}\;s\in I^*\cap (0,\mu(0)),
$$
and conclude that,
$$
\frac{d}{ds}\left[2\gamma_{ E}(K_0)s \frac{dF(s)}{ds}-2\gamma_{ E}(K_0)F(s)+
\frac{1}{2}(F(s))^2 \right]\geq 0,
$$
$\,\mbox{for any}\;s\in I^*\cap (0,\mu(0)).$ Let $P_{+}(s)$ denote the functions in the square brackets.
Since $F$ and $\eta^*$ are both continuous and
locally Lipschitz continuous in $[0,\mu(0)]$ and since, in view of \rife{F2n},
$\frac{dF(s)}{ds}$ is continuous and
locally Lipschitz continuous in $[0,\mu(0)]$ as well,
then we come up with the inequality,
$$
P_{+}(\mu(0))-P_{+}(0)\geq 0.
$$

Therefore we can use \rife{fin-sub1n}, \rife{fin-sub0n} and \rife{F2n} to obtain,
$$
\left[2\gamma_{ E}(K_0)\mu(0)2K_0{\dsp e}^{\eta^*(\mu(0))} -
2\gamma_{ E}(K_0)(2 K_0 M( E))+2(K_0)^2 M^2( E)\right]\geq 0.
$$

Since $\eta^*(\mu(0))=0$, this is equivalent to the following inequality,
$$
2\gamma_{ E}(K_0)\mu(0)-
2\gamma_{ E}(K_0)M( E)+K_0 M^2( E)\geq 0.
$$

So, by using the inequality \rife{Huber-subn} once more and \rife{reg.2-subn} we find,
\beq\label{sharp4}
L^2(\pa E)=\left(\;\int\limits_{\pa E}{\dsp e}^{\frac{u}{2}}ds\right)^2\equiv
\left(\,\int\limits_{\Gamma(0)} d\sigma \right)^2\geq 2\gamma_{ E}(K_0)\mu(0)\geq
\eeq
$$
2\gamma_{ E}(K_0)M( E)-K_0 M^2( E)=
(4\pi-2\mathcal{K}_+(E;K_0)-K_0 M( E))M( E),
$$
which is \eqref{Alex.intron} as claimed.\\

{\bf Step 4.}\\
We will discuss here  the case where the equality holds in \eqref{Alex.intron}.\\
First of all, there is no chance to have the equality in \eqref{Alex.intron} if the strict inequality holds in \eqref{sharp1}.
Therefore, because of \eqref{sharp2}, we see that we must have $\phi=0$ for a.a. $x\in E$, that is, in view
of \rife{13n} and \rife{sharp3}, we also conclude that $u$ must be a solution of \rife{1n} in $E$, and not just a subsolution as in \rife{1nn}.\\
Next we must have the equality sign in the Huber's inequality used in \eqref{Huber-subn} for a.a. $s\in I^*\cap (0,\mu(0))$ and in \rife{sharp4}
for $s=\mu(0)$. Therefore, in view of \eqref{Hubineq-sing} and  \eqref{Hubineq-sing-1}, we conclude that for each
$t\in I\cup\{0\}$, we have,
$$
(a)\qquad \om(t) \mbox{ is simply connected and }\psi(z)=c_t+2\log\left|\Phi_t^{'}(z)(\Phi_t(z))^{-\dsp \al_{\om(t)}} \right|,\;z\in \ov{\om(t)},
$$
where $\al_{\om(t)}=\frac{1}{2\pi}\mathcal{K}_+(\om(t);K_0)$, $c_t\in \R$ and $\Phi_t$ is a conformal map of $\om(t)$ onto the disk of
unit radius $|{\rm w}|=|\Phi_t(z)|<1$ with $\Phi_t(z_t)=0$, for some $z_t \in \om(t)$. Here $\psi$ is the function defined in \eqref{subh-subn}. Since
$\phi$ vanishes, then we have the equality sign in \rife{sharp2} and \rife{sharp1} which therefore do not provide other conditions.
However, in view of the Sard's Lemma for Sobolev functions, we can assume w.l.o.g. that $\om(t)$ is simple for each $t\in I\cup\{0\}$, so that each
$\Phi_t$ can be extended to a univalent and continuous map from $\ov{\om(t)}$ to a closed unit disk, see for example Theorem 2.6 in \cite{pom}.
At this point,
by setting ${\rm w}=\Phi_0(z)$, and in view of $(a)$, we conclude that,
$$
v({\rm w}):=\eta(\Phi^{-1}_0({\rm w})),
$$
is a strong solution of,
$$
-\Delta v=2K_0e^{c_0} |{\rm w}|^{-2 \al}\e{v}\;\; \mbox{in} \;\;\{|{\rm w}|<1\},\qquad v=0  \;\; \mbox{on}\;\; |w|=1,
$$
where $\al=\al_{E}\equiv\al_{\om(0)}$. In particular we have that the level lines of $v$ are concentric circles centred at the origin, that
is, $v$ is radial.
Actually, by using the Brezis-Merle estimates for Liouville type equations (see Remark 5 in \cite{bm}) and standard elliptic theory,
we see that $v$ is analytic far away from the origin and of class $W^{2,q}(B_1)$, for a suitable $q>1$ depending on $\al$.
Thus, by a straightforward evaluation we find that,
$$
v({\rm w})=\log\frac{\tau_0^2}{(1+\frac{K_0 e^{c_0}\tau_0^2}{4(1-\al)^2} |w|^{2(1-\al)})^2},\;|{\rm w}|<1,
$$
for a suitable constant $\tau_0\neq 0$, to be fixed in order to satisfy the Dirichlet boundary condition. As a consequence we find,
$$
\eta(z)=\log\frac{\tau^2 e^{-c_0}}{(1+\frac{ K_0\tau^2}{4(1-\al)^2} |\Phi_0(z)|^{2(1-\al)})^2},\;z\in E,
$$
for some $\tau\neq 0$ and then,
since in particular $e^{\psi(z)}=e^{c_0}\left| \Phi_0^{'}(z)( \Phi_0(z) )^{-\dsp \al}\right|^2$, we see that,
$$
\eta(z)=\log\frac{\tau^2e^{-\psi(z)} \left| \Phi_0^{'}(z)( \Phi_0(z) )^{-\dsp \al}\right|^2}{(1+\frac{K_0\tau^2}{4(1-\al)^2}  |\Phi_0(z)|^{2(1-\al)})^2},\;z\in E.
$$
Since $\eta+\psi=f+u$, then we finally conclude that
$$
\e{f(z)+u(z)}=\frac{\tau^2\left| \Phi_0^{'}(z)( \Phi_0(z) )^{-\dsp \al}\right|^2}{(1+\frac{K_0\tau^2}{4(1-\al)^2}  |\Phi_0(z)|^{2(1-\al)})^2},\;z\in E,
$$
as claimed. Finally, by using the well known fact that the logarithm of the modulus of a non vanishing
holomorphic function is harmonic, we find,

$$
2K e^{f+u}=-\Delta u=
$$
$$
\Delta f -\Delta \log\left(\left| \Phi_0^{'}(z)( \Phi_0(z) )^{-\dsp \al}\right|^2\right)+2\Delta \log\left(1+\frac{K_0\tau^2}{4(1-\al)^2}  |\Phi_0(z)|^{2(1-\al)}\right)=
$$
$$
\Delta f+4\pi \al\dt_{z={z_0}} +2K_0 e^{f+u}= -\omega+4\pi \al\dt_{z=0} +2K_0 e^{f+u},
$$
in the sense of distributions in $E$ and classically in $E\setminus \{0\}$.  Therefore, if $\omega \perp e^{f+u}\mathcal{H}^2$, then this identity
can be satisfied if and only if, 
\beq\label{0308.1}
2K e^{f+u} \equiv 2K_0 e^{f+u}\mbox{ for  a.a. }z\in E, 
\eeq
and $\omega=4\pi \al\dt_{z=z_0}$. In other words
\beq\label{0308.2}
f(z)=h(z)+2\al G(z,z_0)=h(z)-2\al \log |\Phi_0(z)|,
\eeq 
for some $h$ harmonic in $E$. At this point \rife{0308.1} and \rife{0308.2} readily imply that $K\equiv K_0$ for a.a. $z\in E$.
\fineproof

\bigskip
\bigskip

{\bf The proof of  Corollary \ref{co17}.}\\
$(i)$ In this situation we just know 
that $u\in L^{1}_{\rm loc}(\om)$ and $Ke^{f+u}\in L_{\rm loc}^1(\om)$. So we also have $\Delta u\in  L^{1}_{\rm loc}(\om)$ and 
then in particular, by the Green's representation formula, $|\nabla u| \in L^{1}_{\rm loc}(\om)$. 
By Remark 2 in \cite{bm} we find $e^{t|u|}\in L^{1}_{\rm loc}(\om)$ 
for any $t \geq 1$ and letting $\om_0\Subset \om$ be any open, smooth and relatively compact subset, 
we have $u\in L^1(\pa \om_0)$ by standard trace embeddings. Let 
$u=u_1+u_2$, where $u_1$ is the unique weak solution (in the sense of Stampacchia \cite{stam}) of the Dirichlet problem,
$$
\graf{-\Delta u_1=2Ke^{f+u}\quad  \mbox{in}\; \om_0,\\ u_1=0\quad   \mbox{on}\; \om_0,}
$$
and $u_2$ satisfies,
$$
\graf{-\Delta u_2=0 \quad  \mbox{in}\; \om_0,\\ u_2=u\quad   \mbox{on}\; \om_0.}
$$

Then $u_2(x)=-\int\limits_{\pa \om_0} u(y) \frac{\pa G_0}{\pa \nu}(x-y)d\ell_y$, where $G_0$ is the Green's function of $-\Delta$ 
relative to $\om_0$, and since $u\in L^1(\pa \om_0)$, then $u_2\in L^{\ii}_{\rm loc}(\om_0)$. Moreover, 
$u_1\in W_{0}^{1,r}(\om_0)$ for any $r\in (1,2)$ by the results in \cite{stam} and then we find $u\in W_{\rm loc}^{1,r}(\om)$.\\ 

$(ii)$ Let $E$ be any relatively compact and simple subset, 
we can find an open, smooth, simple and relatively compact subset $\om_1$ such that $E\Subset \om_1\Subset \om$. 
Let $K_n\in C^{0}(\ov{\om_1})$ be any sequence satisfying,
\beq\label{0308.3n}
K_n\leq K \mbox{ a.e. in } \om_1\mbox{ and } K_n e^{f+u}\to K e^{f+u},\;\ainf,\mbox{ in } L^{1}(\om_1). 
\eeq
Next, let 
$v_n=v_{n,1}+u_{2}$, where $v_{n,1}$ is the unique weak solution (in the sense of Stampacchia \cite{stam}) of the Dirichlet problem,
$$
\graf{-\Delta v_{n,1}=2 K_n e^{f+u}\quad  \mbox{in}\; \om_1,\\ v_{n,1}=0\quad   \mbox{on}\; \om_1,}
$$
and $u_{2}$ satisfies,
$$
\graf{-\Delta u_2=0 \quad  \mbox{in}\; \om_1,\\ u_2=u\quad   \mbox{on}\; \om_1.}
$$
Obviously, as in $(i)$ we find  $u_2\in L^{\ii}_{\rm loc}(\om_1)$. In particular, by the 
Green's representation formula, it is not difficult to see that,
\beq\label{0408.1n}
v_n\leq u \mbox{ a.e. in } \om_1.
\eeq
Let us observe that, by Theorem \ref{thm2}, 
$e^{f+u}=e^{\rho}\in L^{p_0,{\rm loc}}_{\rm loc}(\om\setminus S_{2\pi})\cap L^{q_0}_{\rm loc}(\om)$ for some $p_0>2$ and $q_0>1$, 
whence by standard elliptic estimates and the Sobolev embedding 
we find $v_n\in W^{2,q_0}(\om_1)\cap C^{0}(\ov{\om_1})$. 
By using \rife{0308.3n} with well known results in \cite{stam}, we conclude that $v_n\to u$ in $W^{1,r}_{\rm loc}(\om_1)$, for any $r\in (1,2)$. 
At this point we observe that $v_n$ is a solution of, 
$$
-\Delta v_{n}=2\widehat{K_n} e^{f}e^{v_n}  \mbox{ in }\om_1,
$$
where, 
$$
\widehat{K_n}=K_n e^{u-v_n}\mbox{ satisfies } \sup\limits_{\ov{\om_1}}|\widehat{K_n}|\leq C_n e^{u}.
$$
By $(i)$ we have $\widehat{K_n}\in L^{t}(\om_1)$ for any $t\geq 1$. On the other side, by Proposition 
\ref{prexp}, we also find that $e^f \in L^{s,{\rm loc}}_{\rm loc}(\om\setminus S_{2\pi})\cap L^{m}_{\rm loc}(\om)$ for some $s>2$ and $m>1$.
Therefore we can apply Theorem \ref{thm1}$(i)$ on $\om_1$ with $K=\widehat{K_n}$ and $u=v_n$, to conclude that, 
\beq\label{Alexappn}
\left(\,\int\limits_{\pa E} e^{\frac{f+v_n}{2}}d\ell\right)^2\geq 
\left(4\pi-2\mathcal{K}_{+,n}(E;K_0)-K_0\int\limits_{ E }e^{f+v_n}\right)\int\limits_{ E }e^{f+v_n},
\eeq
where, 

$$
\mathcal{K}_{+,n}(E;K_0)=k_{s,+}(E)+\int\limits_{E}[\widehat{K_n}-K_0]^+{\dsp e}^{f+v_n}dx.
$$

Since  $v_n\to u$ in $W_{\rm loc}^{1,r}(\om_1)$ and in view of \rife{0308.3n}, 
along a subsubsequence (which we will not relabel) we have $v_n\to u$ a.e. in $\om_1$ and $K_n e^{f+u}\to K e^{f+u},\;\ainf,\mbox{ a.e. in }\om_1$.
Then, by \rife{0408.1n} and the dominated convergence theorem we conclude that, 
$$
\int\limits_{ E }e^{f+v_n}\to \int\limits_{ E }e^{f+u},\ainf,
$$
$$
\mathcal{K}_{+,n}(E;K_0)\to\mathcal{K}_{+}(E;K_0),\ainf,
$$
and, 
$$
\int\limits_{\pa E} e^{\frac{f+v_n}{2}}d\ell\to \int\limits_{\pa E} e^{\frac{f+u}{2}}d\ell,
$$
where for the second limit we observe that, 
$$
[\widehat{K_n}-K_0]^+{\dsp e}^{f+v_n}= [K_n e^{u-v_n}-K_0]^+{\dsp e}^{f+v_n}\leq 
[K e^{u-v_n}]^+{\dsp e}^{f+v_n}= [K ]^+{\dsp e}^{f+u}.
$$
It is understood that the last limit holds true whenever $\int\limits_{\pa E} e^{\frac{f+u}{2}}d\ell$ is finite, otherwise  
\rife{Alex.intron} is trivially satisfied since $M(E)<+\ii$.
Therefore, in the limit $n\to +\ii$, along the given subsequence we recover \rife{Alex.intron}, as claimed. \fineproof

\bigskip

\section{Application to the Alexandrov's isoperimetric inequality.}\label{sec4}
The notion of Surface of Bounded Curvature ({\it SBC} for short) was introduced by A.D. Alexandrov \cite{Ale}, as a
model to describe surfaces with a wide variety of singularities.  A detailed discussion of this subtle subject is behind the scope
of our work, and we refer the reader to \cite{Ale2} and \cite{Res2} for a complete account about the subject, and to \cite{Troy1} for a shorter
exposition of some of the main results.
Here we will just use an equivalent local description of
these objects.\\
Indeed, according to a series of results due to Huber and Reshetnyak, see  \cite{Res2},
an {\it SBC} without boundary can be equivalently defined as a Riemann surface $\mathcal{M}$ equipped with a metric $\mathfrak{g}$, which admits an
atlas of local charts $\mathcal{U}=\{U_j, \Phi_j \}_{j\in J}$, such that each $\Phi_j$ is an isometry of $U_j$ on $\om_j=\Phi_j(U_j)$, with
$\om_j\subset \RR^2(\simeq \CC)$, a smooth, open and bounded set, such that $\mathfrak{g}$ in local coordinates
takes the form of a quadratic differential, $\Phi_j^{\mbox{\#}}(\mathfrak{g})= e^{\rho_j(z)}|dz|^2$,  $z=x+iy\in \CC$. 
Here  $\mbox{\#}$ denotes the standard pull-back, 
$|dz|^2$ is the Euclidean metric and $\rho\equiv \rho_j=\rho_+-\rho_-$,
where $\rho_\pm$ are two superharmonic functions defined by,
\beq\label{frepp}
\rho_{\pm}(z)=h^0_{\pm}(z)+\int\limits_{\om_j} \Gamma (z,y) d\omega^{0}_{\pm}(y),\quad \Gamma(z,y)=\frac{1}{2\pi}\log\left(\frac{1}{|z-y|}\right),
\eeq
with $h^0_{\pm}$ harmonic in $\om_j$. Here $\omega^0_\pm$ are the mutually orthogonal non negative measures
defined by the Jordan decomposition of a measure of bounded total variation on $\om_j$,  $\omega^0=\omega^{0}_{+}-\omega^{0}_{-}$ .
Any such system of coordinates is said to be isothermal and any metric taking the form $e^{\rho(z)}|dz|^2$ with $\rho$ as in \rife{frepp}
is said to be subharmonic. Among other things,
the definition is completed by the transitions rules between charts of functions and holomorphic forms, thus including the metric, see \cite{Res2}
for further details.\\

This is why we will focus our attention on the local model of an {\it SBC}.
\bdf {\it An Abstract Surface of Bounded Curvature ({\it ASBC} for short) is a pair
$\mathcal{S}=\left\{\om, e^{\rho(z)}|dz|^2\right\}$, where $\om\subset \RR^2$ is open, smooth and bounded and
$\rho=\rho_+-\rho_-$, with $\rho_{\pm}$ as defined in } \rife{frepp}.
\edf

So, if $\mathcal{S}=\left\{\om, e^{\rho(z)}|dz|^2\right\}$ is an {\it ASBC},
according to Reshetnyak (see \cite{Res2} Theorem 8.1.7), the \uv{total curvature} $\mathcal{K}$, is the measure of finite total variation defined as follows,
\bdf\label{300616.2}{\it Let $\mathcal{S}=\left\{\om, e^{\rho(z)}|dz|^2\right\}$ be an {\it ASBC}.
The total curvature $\mathcal{K}(E)$ of a Borel set $E\subseteq \om$ is defined by:
$$
2\mathcal{K}(E):=\omega^0(E)=\omega^0_+(E)-\omega^0_-(E).
$$}
\edf

\bigskip

\brm{\it We remark that, with this definition, the total curvature is well defined and finite for any Borel set $E\subseteq \om$.
Nevertheless, if for some $z_0\in\om$ it holds $\omega^{0}_{+}(z_0)\geq 4\pi$,  then the lengths and areas of sets containing $z_0$, as
defined via the metric $g=e^{\rho(z)}|dz|^2$ {\rm(}see {\rm \rife{elle}}, {\rm \rife{emme}} below{\rm ) } are not well defined in general.
Any point $z_0\in\om$ which satisfies  $\omega^{0}_{+}(z_0)\geq 4\pi$ is said
to be a \uv{cusp}.}
\erm

From now on we will assume that $\mathcal{S}=\left\{\om, e^{\rho(z)}|dz|^2\right\}$ is an {\it ASBC}
\uv{with no cusps}, that is, we assume that,
\beq\label{nocusps}
\;\fo\,z\in\om,\;\omega^{0}_{+}(z)< 4\pi.
\eeq

Let $S_{2\pi}=\{x\in \om\;:\;\omega^{0}_{+}(z)\geq 2\pi\}$. We have the following seemingly new result about the structure of
subharmonic metrics with no cusps. Interestingly enough it is sharp, see Example 1 below for further details. 
The proof is based on various results and arguments in \cite{bm} about the regularity properties of Liouville-type equations.
Here $\mathcal{H}^\gamma$, with $\gamma>0$, denotes the $\gamma$-dimensional Hausdorff measure.
\bte\label{thm2}
Let  $\mathcal{S}=\left\{\om, e^{\rho(z)}|dz|^2\right\}$ be an {\it ASBC} with no cusps.\\
Then
$e^{\rho}\in L^{p_0,{\sscp \rm loc}}_{\rm loc}(\om \setminus S_{2\pi})\cap L^{q_0}_{\rm loc}\left(\om\right)$ for some $p_0>2$ and some $q_0>1$.
Moreover, there exists $K\in L^1_{\rm loc}\left(\om;e^{\rho}\mathcal{H}^2\right)$ 
and a Radon measure $k_s$ on $\om$, satisfying $k_s\perp e^{\rho}\mathcal{H}^2$, such that, letting
$k_s=k_{s,+}-k_{s,-}$ be the Jordan decomposition of $k_s$, then
$\rho$ can be decomposed as $\rho=u+f$, where $f=f_+-f_-$, with $f_\pm$ satisfying {\rm \rife{frep}} with  $\omega_\pm=2k_{s,\pm}$ and
$h_\pm$ suitable harmonic functions and where $u\in L^{1}_{\rm loc}(\om)$ is a solution of,
\beq\label{kreg}
-\Delta u=2 K e^{f+u}\;\;\mbox{ in }\;\;\om,
\eeq
in the sense of distributions. In particular, either,\\
$(i)$ $K\in L^s_{\rm loc}\left(\om\right)$, for some $s>\frac{q_0}{q_0-1}$
and then $u$ is a strong  solution of {\rm \rife{kreg}} which satisfies
$u\in L^{\ii}_{\rm loc}(\om)\cap W_{\rm loc}^{2,r}(\om) $, $\fo\;1\leq  r\leq \frac{s q_0}{s+q_0}$, or,\\
$(ii)$ $u\in W_{\rm loc}^{1,r}(\om)$ for any $r\in (1,2)$ and $e^{t|u|}\in L^{1}_{\rm loc}(\om)$ for any $t \geq 1$.\\
In both cases,
$$
e^{\rho(z)}|dz|^2\equiv e^{u(z)+f(z)}|dz|^2,\;z\in\om,\quad K e^{f+u}\in L^{1}_{\rm loc}(\om),
$$
and
\beq\label{curv0}
\mathcal{K}(E)=\int\limits_E K e^{f+u}+ k_s(E),
\eeq
for any relatively compact Borel set $E\Subset \om$. Moreover, if $\rho=u+f$ for a pair $\{u,f\}$ as above, then, for any fixed $h$ harmonic in $\om$,
the pair $\{u_h,f_h\}:=\{u-h,f+h\}$ satisfies the same properties with $\rho=u_h+f_h$.
\ete

\proof 
Let $H(z,y)=G(z.y)-\Gamma(z,y)$ be the regular part of the Green's function on $\om$. Then 
${m}_{\pm}(z)=\ino H(z,y)d\omega^{0}_{\pm}(y)$ are harmonic in $\om$, and $\rho+m_+-m_-$ takes the form $\rho_+-\rho_-$ for a suitable pair 
$\rho_{\pm}$ satisfying \rife{frep}.
Therefore, by Proposition \ref{prexp}, we find
$e^{\rho}\in L^{p_0,{\sscp \rm loc}}_{\rm loc}(\om \setminus S_{2\pi})\cap L^{q_0}_{\rm loc}\left(\om\right)$ for some $p_0>2$
and some $q_0>1$. Then $e^{\rho} \mathcal{H}^2$ is a Radon measure on $\om$, and so it is well defined
the Lebesgue decomposition of
$\mathcal{K}$ with respect to $e^{\rho} \mathcal{H}^2$,
\beq\label{curv}
\mathcal{K}=K e^{\rho}\mathcal{H}^2+k_s,\quad K\in L_{\rm loc}^1\left(\om;e^{\rho}\mathcal{H}^2\right),\quad k_s\perp e^{\rho} \mathcal{H}^2,
\eeq
where $k_s$ is a Radon measure on $\om$.
We first observe that, since $\rho\in L_{\rm loc}^1(\om)$, then
$-\Delta \rho =\omega_+^0-\omega_-^0$ holds in the sense of distributions in $\om$, whence,
by \rife{curv} and the definition of $\mathcal{K}$, we see that the following equality,
$$
-\Delta \rho = 2Ke^{\rho} +2 k_s,
$$
holds as well, in the sense of distributions in $\om$. Let $f=f_{+}-f_{-}$ be defined by \rife{frep} with $\omega_\pm=2k_{s,\pm}$, $h_\pm=0$,
and let us set,
$$
u:=\rho-f.
$$
Clearly $u\in L^{1}_{\rm loc}(\om)$, and since $-\Delta f=2k_s$ in the sense of distributions, then we deduce that,
$$
-\Delta u = 2Ke^{f+u} +2 k_s +\Delta f=2Ke^{f+u},
$$
that is, $u$ satisfies \eqref{kreg} in the sense of distributions in $\om$.\\
At this point, the fact that $\mathcal{K}(E)$ takes the form \rife{curv0} is a straightforward consequence of the fact that $k_s\perp e^{\rho} \mathcal{H}^2$.
Moreover we observe that, if $K$ satisfies the assumption in $(i)$, then all the assumptions of Lemma \ref{lem0308} $(a)$ are satisfied 
and then the conclusion readily follows.\\

So we are left with the case where $K$ does not satisfy the assumption in $(i)$, that is, 
we just know that $K\in L^1(\om; e^{f+u}\mathcal{H}^2)$ where 
$u\in L^1\left(\om\right)$ is a solution of \eqref{kreg} in the sense of distributions. Therefore all the assumptions of Corollary \ref{co17}$(i)$ are 
satisfied and then the desired conclusion follows.\\
Finally it is obvious that the representation $\rho=u+f$ with all the properties established above still holds for $\{u_h,f_h\}$ where $h$ is 
an arbitrary harmonic function in $\om$.
\fineproof

\bigskip
\bigskip

Let $E\Subset \om$ be any regular and relatively compact subset and suppose that \eqref{nocusps}
holds. Then we define
the length of $\pa E$,
\beq\label{elle}
L(\pa E)=\int\limits_{\pa E} e^{\frac{f+u}{2}}d\ell,
\eeq
and the area of $E$,
\beq\label{emme}
M( E )=\int\limits_{ E }e^{f+u}dx.
\eeq

\bdf {\it
For $K_0\in \R$ and for any and relatively compact Borel set  $E\Subset \om$, we
define the positive variation of the total curvature of $E$ with respect to $K_0$,
\beq\label{gammadefnn}
\mathcal{K}_{+}(E;K_0)=\sup\limits_{U\subseteq E} \left \{\mathcal{K}(U)-K_0\int\limits_{U}{\dsp e}^{f+u}dx\right \},
\eeq
where the supremum is taken over all Borel sets $U\subseteq E$.}
\edf
Because of \rife{curv0}, and since $k_{s,+}\perp{\dsp e}^{f+u}\mathcal{H}^2$, then
$\mathcal{K}_+ (E;K_0)$ takes the form,
$$
\mathcal{K}_+(E;K_0)=k_{s,+}(E)+\int\limits_{E}[K-K_0]^+{\dsp e}^{f+u}dx.
$$

\bdf\label{defcone}{\it
For fixed $\al>-1$ and $K_0>0$, a spherical $\{K_0,\al\}$-cone is the {\it ASBC} defined by 
$\left\{B_1, |{\rm w}|^{-2\al}e^{v({\rm w})}|d{\rm w}|^2\right\}$
where $B_1=\{{\rm w}\in \ci\,:\,|{\rm w}|<1\}$ and,
$$
e^{v({\rm w})}=\frac{\tau_0^2}{\left(1+\frac{K_0 \tau_0^2}{4(1-\al)^2} |{\rm w}|^{2(1-\al)}\right)^2},\;|{\rm w}|<1,
$$
for some $\tau_0\neq 0$.}
\edf

It is worth to remark that the function $v$ in Definition \ref{defcone} is of class $L^{\ii}(B_1)\cap W^{2,p}_{\rm loc}(B_1\setminus \{0\}\})\cap W^{2,q}(B_1)$
for any $p>2$ and for any $q<\frac{1}{|\al|}$ and it is a strong solution of $-\Delta v = 2K_0 |{\rm w}|^{-2\al}e^{v}$ in $B_1$.

\bigskip

In view of Theorem \ref{thm1}, Corollary \ref{co17} and Theorem \ref{thm2}, and in the same spirit of \cite{Band0}, for $K_0\geq 0$ we obtain a new proof of the Alexandrov \cite{Ale} isoperimetric inequality 
on an ASBC. At least to our 
knowledge the characterization of the equality sign in this weak framework is new. 
\bte\label{Alexandrovn}
Let $\mathcal{S}=\left\{\om, e^{\rho}|dz|^2\right\}$ be an {\it ASBC}
with no cusps and fix $K_0\geq 0$. Then the curvature takes the form {\rm \rife{curv0}} for some $u,f,K,k_s$ as in Theorem \ref{thm2} and
for any simple and relatively compact subset $E\Subset \om$, it holds,
\beq\label{Alex.intro}
L^2(\pa E)\geq \left(4\pi-2\mathcal{K}_+(E;K_0)-K_0M(E)\right)M(E).
\eeq

In particular, if $K$ satisfies the assumption of Theorem \ref{thm2}$(i)$ and also $K\in  L^{n,{\sscp \rm loc}}_{\rm loc}(\om\setminus S_{2\pi}) $, 
for some $n>\frac{2 p_0}{p_0-2}$, then the equality in {\rm \rife{Alex.intro}} holds if and only if:\\ 
- $\left\{E, e^{\rho}|dz|^2\right\}$ is isometric to a spherical $\{K_0,\al\}$-cone with $\al=\frac{1}{2\pi}k_{s,+}(E)$;\\
- $\rho=u+f$ and $e^{f+u}$ takes the form {\rm \rife{sharp-last}}, where $u$ is a solution of {\rm \rife{kreg}} with $K\equiv K_0$ for $a.a.\;z\in E$
and $k_s=2\pi \al \dt_{z=z_0}$, for some $z_0\in E$, that is, $f(z)=h(z)+2\al G(z,z_0)=h(z)-2\al \log|\Phi_0(z) |$, for some function
$h$ harmonic in $E$.
\ete
\proof Since $\mathcal{S}$ is an {\it ASBC} with no cusps, then, by Theorem \ref{thm2},
the curvature takes the form \rife{curv0} where $u\in L^{1}_{\rm loc}(\om)$ is a solution of \rife{kreg} in the sense of distributions, 
$f$ takes the form \rife{frep} with $h_{\pm}$ harmonic and $\omega_{\pm}=2k_{s,\pm}$ and $Ke^{u+f}\in L^1(\om)$. 
If $K$ satisfies the assumption in Theorem \ref{thm2}$(i)$ and also $K\in  L^{n,{\sscp \rm loc}}_{\rm loc}(\om\setminus S_{2\pi}) $, 
for some $n>\frac{2 p_0}{p_0-2}$, then all the hypothesis
of Theorem \ref{thm1}$(i)$ are satisfied as well.
As a consequence, the inequality \rife{Alex.intro} holds and the equality sign is attained if and only if
\rife{sharp-last} holds, that is,
$$
e^{f(z)+u(z)}|dz|^2=|\sg \Phi_0(z)|^{-2\al}e^{v(\sg\Phi_0(z))}|d(\sg\Phi_0(z))|^2=|{\rm w}|^{-2\al}e^{v({\rm w})}|d{\rm w}|^2,\;\sg=\sqrt[1-\al]{\tau},
$$
for any $B_1\ni {\rm w}=\Phi_0(z),\;z\in E$, as claimed. In particular,   since $\omega=2k_s \perp e^{u + f}\mathcal{H}^2$ by construction, then 
Theorem \ref{thm1}$(iii)$ can be applied as well. This observation completes the discussion of the equality case.\\ 
Clearly, to conclude the proof, it is enough to show that \eqref{Alex.intro} holds in case $(ii)$ of Theorem \ref{thm2} is satisfied. 
However this is just the content of Corollary \ref{co17}$(ii)$ which immediately yields the desired conclusion.
\fineproof

\bigskip
\bigskip

\section{Examples.}\label{sec5}

We recall that a point $P$ on an {\it SBC} is said to be a conical singularity of order $\al>-1$ if in an isothermal chart
$\{\om,z\}$ such that $z(P)=0$, the metric takes the form $e^{\rho(z)}|dz|^2= |z|^{2\al} e^{u(z)}|dz|^2$, where
$u\in C^{0}(\om)\cap C^{2}(\om\setminus \{0\})$.\\
In this section $\dt_p$ denotes the Dirac delta with pole at $p\in \R^2$.

\bigskip

{\bf Example 1.}\\ 
We use Example 1 in \cite{bm} to construct an ASBC 
$\left\{B_1, e^{\rho}|dz|^2\right\}$ such that $\{u,f,K,k_s\}$ as obtained in Theorem \ref{thm2} have the following properties:\\
- either $e^{\rho}\in L^\ii(B_1)$ or $e^{\rho}\in L^q(B_1)$, for any $q\geq 1$;\\
- $K\in L^{1}(e^\rho\mathcal{H}^2,B_1)\cap L^{1}(B_1)$ but there is no $s>1$ such that $K\in L^{s}(B_1)$;\\
- $u$ is not locally bounded;\\
- $u$ has all the properties claimed in Theorem \ref{thm2}$(ii)$.\\
Let $0\neq a<1$,  and for $z\in B_1\setminus\{0\}$ let us set,
$u(z)=-a \log\left(\log\left(\frac{e}{|z|}\right)\right)$ and, 
$$
K(z)=-\frac{a}{2} |z|^{-2} \left(\log\left(\frac{e}{|z|}\right)\right)^{-(2-a)}.
$$
The superharmonic function $\rho(z)=\int\limits_{B_1}G(z,y)d\omega^0(y)$, where $\omega^{0}(y)=2K(y) e^{u(y)}d\mathcal{H}^2$, 
takes the form $\rho=\rho_+-\rho_-$ as in \rife{frepp} with $h^0_{\pm}=0$, and 
$\omega^0_-=0$ and $\omega^{0}_+(y)=2K(y) e^{u(y)}d\mathcal{H}^2$ if $a<0$, while 
$\omega^{0}_-(y)=2K(y) e^{u(y)}d\mathcal{H}^2$ and $\omega^0_+=0$ if $a\in (0,1)$.
Since $Ke^{u}\in L^{1}(B_1)$, then $\omega^{0}<<e^{u}\mathcal{H}^2$ and so we find $\{u,f,K,k_s\}$ as claimed in Theorem \ref{thm2} by setting 
$f=0$, $k_s=0$,  $K\equiv K$ and $u\equiv u$. In fact we see that $u$ is a solution of, 
$$
\graf{-\Delta u=2K e^{u}\quad  \mbox{in}\; B_1,\\ u=0\quad   \mbox{on}\; B_1,}
$$
that is, in particular $u\equiv \rho$, and so we find,
$$
e^{\rho(z)}= \left(\log\left(\frac{e}{|z|}\right)\right)^{-a},\; z\in B_1.
$$

If $a\in(0,1)$, then $e^\rho\in L^{\ii}(B_1)$,  $K\in L^{1}(e^\rho\mathcal{H}^2,B_1)\cap L^{1}(B_1)$ but $u(z)\to -\ii$ as $z\to 0$. 
If $a<0$, then $e^\rho\in L^{q}(B_1)$ for any $q\geq 1$, $K\in L^{1}(e^\rho\mathcal{H}^2,B_1)\cap L^{1}(B_1)$ but $u(z)\to +\ii$ as $z\to 0$. 
In both cases, there is no $s>1$ such that $K\in L^{s}(B_1)$, so there is no chance that $K$ satisfies the assumption of Theorem \ref{thm2}$(i)$. 
On the other side, in both case it is easy to check that $u$ has all the properties claimed in Theorem \ref{thm2}$(ii)$.
Clearly Theorem \ref{Alexandrovn} applies and then \rife{Alex.intro} holds on $\left\{B_1, e^{\rho}|dz|^2\right\}$.

\bigskip
\bigskip

{\bf Example 2.}\\
Let $\mathbb{S}^{2}_{\al_1,\al_2}$ be the {\it SBC} defined by the isothermal charts $\{\om_i,\varphi_i\}_{i=1,2}$ and the local metrics
$\{g_i\}_{i=1,2}$ constructed as follows. For $r_0 \geq 4$ and $-1<\al_1\leq \al_2\leq 0$, we define,
$$
\om_2=\{z\in \ci\,:\,|z|<r_0\},\quad  \varphi_2=z, \quad g_2=e^{\rho}|dz|^2,
$$

$$
\om_1=\left\{z\in \ci\cup \{\ii\}\,:\,|z|>\frac{1}{r_0}\right \},\quad  \varphi_1=\frac{1}{z}, \quad g_1=\varphi_1^{\mbox{\#}}(g_2),
$$
where,
\beq\label{form}
\rho(z)=
\graf{
\log{\left(\dsp \frac{ 4 (1+\al_2)^2  |z|^{2\al_2}}{ \left( 1 + |z|^{2(1+\al_2)}\right)^2 } \right)},
&\,\,|z|<1,\\\\
\log{\left(\dsp \frac{ 4 (1+\al_2)^2 |z|^{2\al_1}}
{ \left( 1 + |z|^{2(1+\al_1)}\right)^2 }\right)},&\,\,|z|\in\,[1,+\infty).
}
\eeq

This is a compact surface without boundary,
homeomorphic to the two sphere, with two conical singularities, $z=\ii$ of order $\al_1$ and $z=0$ of order $\al_2$.
For $\al_1=\al_2<0$ we are reduced to the classical "american football" \cite{Troy0}, with constant Gaussian curvature $K\equiv 1$.
Instead, if $\al_1<\al_2\leq 0$, we have the glueing of two caps of american footballs
with gaussian curvatures $1$ and $\frac{(1+\al_1)^2}{(1+\al_2)^2}$ respectively, with different conical singularities,
see \cite{barjga} and \cite{Fang} for more details about this singular surface.\\
We consider a decomposition in the $\{\om_2,\varphi_2\}$ chart, as claimed in Theorem \ref{thm2}, of the form
$
\rho(z)=f(z)+u(z),
$
where,
$$
u(z)=
\graf{
\log{\left(\dsp \frac{ 4 (1+\al_{2})^2 }{ \left( 1 + |z|^{2(1+\al_2)}\right)^2 } \right)},
&\,\,|z|<1,\\\\
\log{\left(\dsp \frac{ 4 (1+\al_{2})^2 |z|^{2(\al_1-\al_2)}}
{ \left( 1 + |z|^{2(1+\al_1)}\right)^2 }\right)},&\,\,|z|\in\,[1,+\infty),
}
$$
and
$$
f(z)=f(z;\al_2)=2\al_2\log{|z|},\; |z|\in (0,+\ii).
$$

Clearly we have $u\in W_{\rm loc}^{2,k}(\R^2\setminus \{0\})\cap W_{\rm loc}^{2,r}(\R^2)$, for any $k>2$ and $1<r<\frac{2}{|\al_1|}$,
which is also a strong solution of $-\Delta u = 2K |z|^{2\al_2}e^{u}\;\mbox{in}\;\R^2$, with,
$$
K(z)=\graf{1, & |z|\in [0,1),\\ \frac{(1+\al_1)^2}{(1+\al_2)^2}, & |z|\in(1,+\ii).}
$$

So $K\in L^{\ii}(\R^2)$ and putting,
$$
k_{s,2}=2\pi |\al_2| \delta_{z=0},
$$
we find,
$$
\mathcal{K}(E)=\int\limits_E K e^{ f+u}d\mathcal{H}^2 + k_{s,2}(E),\quad E\Subset \{  |z| < r_0\},
$$
which is the total curvature of a relatively compact Borel set $E$ in the $\{\om_2,\varphi_2\}$ chart. 
For a generic Borel set $E_0 \subseteq \ci\cup \{\ii\}$, we can consider the analogue decomposition for $g_1$ which takes the form 
$g_1=e^{\rho_1}|d{\rm w}|^2$, with $\rho_1=f_1+u_1$, where 

$$
u_1({\rm w})=
\graf{
\log{\left(\dsp \frac{ 4 (1+\al_{2})^2 }{ \left( 1 + |{\rm w}|^{2(1+\al_1)}\right)^2 } \right)},
&\,\,|{\rm w}|<1,\\\\
\log{\left(\dsp \frac{ 4 (1+\al_{2})^2 |{\rm w}|^{2(\al_2-\al_1)}}
{ \left( 1 + |{\rm w}|^{2(1+\al_2)}\right)^2 }\right)},&\,\,|{\rm w}|\in\,[1,+\infty),
}
$$
$f_1(\cdot)=f (\,\cdot\,, \al_1)$, and 
eventually find the total curvature of any Borel set $E_0 \subseteq \ci \cup \{\ii\}$,
\beq\label{last}
\mathcal{K}(E_0)=\int\limits_{E_{0,2}} K e^{f+u}d\mathcal{H}^2 + k_{s,2}(E_{0,2}) +
\int\limits_{\varphi_1(E_{0,1})} K_1 e^{f_1 +u_1}d\mathcal{H}^2 + k_{s,1}(\varphi_1(E_{0,1})), 
\eeq
where $E_{0,2}=E_0 \cap \{|z|< r_0 \}$, $E_{0,1}=E_0 \cap \{|z|\geq r_0 \}$, $K_1=K\circ \varphi_1$, 
and, 
$$
k_{s,1}=2\pi |\al_1| \delta_{{\rm w}=0}.
$$

Next, to simplify the notations let us set,
$$
\sg_{1,2}=\frac{(1+\al_1)^2}{(1+\al_2)^2}\leq 1.
$$

It is easy to check that the area of $\mathbb{S}^{2}_{\al_1,\al_2}$ is $2\pi(1+\al_2)+\frac{1}{\sg_{1,2}} 2\pi(1+\al_1)$ while, by using \rife{last}, we see that 
the total curvature of $\mathbb{S}^{2}_{\al_1,\al_2}$ is $4\pi$, in agreement with the fact that, as well known \cite{Ale2}, 
the Gauss-Bonnet formula holds even in this singular context. Please observe that this is just an equivalent formulation of
the singular Gauss-Bonnet formula, see \cite{Troy}, which asserts that the global integral of the absolutely continuous part of the Gaussian curvature equals the 
singular Euler characteristic, yielding in this particular case the well known identity,
$$
\int\limits_{B_1} K e^{f+u}d\mathcal{H}^2+\int\limits_{\varphi_1((B_1)^c)} K_1 e^{f_1 +u_1}d\mathcal{H}^2=2\pi(2+\al_1+\al_2).
$$

If $E$ is a simple set surrounding the origin, then we can always take $r_0$ large enough to guarantee that $E\Subset \{  |z| < r_0\}$ so that 
the inequality \rife{Alex.intro} takes the form,
$$
L^2(\pa E)\geq \left(4\pi(1+\al_2) -2[1-K_0]^+M(E\cap B_1)-2[\sg_{1,2}-K_0]^+M(E\cap (B_1)^{c})-K_0M(E)\right)M(E).
$$
In particular, if $K$ is not constant in $E$, then the inequality is always strict and if we choose $K_0= 1$, then
it reduces to the well known Bol's \cite{Bol} inequality,
$$
L^2(\pa E)\geq \left(4\pi(1+\al_2) -M(E)\right)M(E).
$$
If $E=B_R$ with $R\leq 1$, then $K\equiv 1$ in $E$ and since,
$$
L^2(\pa B_R)=\left(\int_0^{2\pi}\frac{2 (1+\al_2) R^{\al_2}}{1 + R^{2(1+\al_2)}}d\ell \right)^2=
\frac{16\pi^2 (1+\al_2)^2  R^{2\al_2}}{( 1 + R^{2(1+\al_2)})^2},
$$
and,
$$
M(B_R)=\int\limits_{B_R}\dsp \frac{ 4 (1+\al_2)^2  |x|^{2\al_2}}{ \left( 1 + |x|^{2(1+\al_2)}\right)^2 } dx=
\frac{4\pi (1+\al_2) R^{2\al_2}}{ 1 + R^{2(1+\al_2)}},
$$
then we find the equality in \rife{Alex.intro} with $K_0=1$,
$$
L^2(\pa B_R)=\frac{16\pi^2 (1+\al_2)^2  R^{2\al_2}}{( 1 + R^{2(1+\al_2)})^2}=
$$
$$
\left(4\pi(1+\al_2) -\frac{4\pi (1+\al_2)  R^{2\al_2}}{ 1 + R^{2(1+\al_2)}}\right)\frac{4\pi (1+\al_2)  R^{2\al_2}}{ 1 + R^{2(1+\al_2)}}=
\left(4\pi(1+\al_2) -M(B_R)\right)M(B_R).
$$

\bigskip
\bigskip

{\bf Example 3.}\\
This example illustrates the failure of Theorem \ref{thm2} on a surface homeomorphic to the two-sphere with a cusp and in the same time
the kind of singularity which yields a curvature function $K$ which is unbounded but in $L^r(E)$ for some $r>1$.\\

Let us consider the same charts $\{\om_i,\varphi_i\}_{i=1,2}$ as in Example 2, where this time the metric $g_2(z)=e^{\rho(z)}|dz|^2$ is defined as follows,
\beq\label{form2}
\rho(z)=
\graf{
\log{\left(\dsp \frac{ 2 }{ \left( 2 - |z|^{\frac12}\right)^2 } \right)},
&\,\,|z|<1,\\\\
\log{\left(\dsp \frac{8 |z|^{\frac32}}{ \left( 1 + |z|^{\frac12}\right)^2 }\right)},&\,\,|z|\in\,[1,+\infty).
}
\eeq

We consider a decomposition as claimed in Theorem \ref{thm2} in the $\{\om_2,\varphi_2\}$ chart,
$
\rho(z)=f(z)+u(z),
$
where we choose $f=0$ so that $u=\rho$, which satisfies
$u\in W_{\rm loc}^{2,k,{\sscp \rm loc}}(\R^2\setminus \{0\})\cap W_{\rm loc}^{2,r}(\R^2)$, for any $k>2$ and $1<r<\frac{4}{3}$,
and is a strong solution of $-\Delta u = 2Ke^{f+u}\;\mbox{in}\;\R^2$, where,
$$
K(z)=\graf{-\frac14 \frac{1}{|z|^{\frac32}}, & |z|\in [0,1),\\ \frac{1}{32} \frac{1}{|z|^3}, & |z|\in(1,+\ii).}
$$

The total curvature of a relatively compact Borel set $E$ in the $\{\om_2,\varphi_2\}$ chart takes the form,
$$
\mathcal{K}(E)=\int\limits_E Ke^{u}dx,\quad E\Subset  \{  |z| < r_0\},
$$
with $K\in L^r(E)\cap L^{\ii}_{\rm loc}(E\setminus \{0\})$ for any $1<r<\frac{4}{3}$.\\
On the other side, let us check whether or not the assumption \rife{nocusps} is satisfied on a generic relatively compact Borel set in the $\{\om_1,\varphi_1\}$ chart.
The metric takes the form,
$$
g_1({\rm w})=\varphi_1^{\mbox{\#}}(g_2)=e^{\rho_1({\rm w})}|d{\rm w}|^2,
$$
where,
\beq\label{form3}
\rho_1({\rm w})=
\graf{
\log{\left(\dsp \frac{8 |{\rm w}|^{-\frac92}}{ \left( 1 + |{\rm w}|^{\frac12}\right)^2 }\right)},&\,\,|{\rm w}|\in\,[0,1],\\\\
\log{\left(\dsp \frac{ 2 |{\rm w}|^{-3}}{ \left( 2 |{\rm w}|^{\frac12} - 1\right)^2 } \right)},
&\,\,|{\rm w}|\in\,(1,+\infty).
}
\eeq

Therefore, it is readily seen that $\rho_1$ takes the form \rife{frepp} with $\omega_+^{0}(0)=\frac{9\pi}{2}>4\pi$ which violates \rife{nocusps}.
This singular surface is still homeomorphic to the two sphere, but it has a cusp at $z=\ii$.
As a consequence, while the curvature is always well defined in the sense of measures, the area of a compact Borel set in the
$\{\om_1,\varphi_1\}$ chart is not, since $e^{\rho_1} $ is not an $L^1_{\rm loc}(\R^2)$ function. In particular,
there is  no chance to use
the argument in the proof of Theorem \ref{thm2}, which should be based on the Lebesgue decomposition of $\mathcal{K}=\omega^0$ with respect to
$e^{\rho_1}\mathcal{H}^2$, since the latter is not even a Radon measure in this case. It is worth to mention that, nevertheless, the product
$(K\circ \varphi_1) e^{\rho_1}$ is an $L^1_{\rm loc}(\R^2)$ function which could be used in principle as the density of the total curvature. 
On the other hand, the right hand side of the Alexandrov's isoperimetric inequality \rife{Alex.intro} is not well defined in general.\\
However Theorem \ref{thm2} and
Theorem \ref{Alexandrovn} can be applied in the  $\{\om_2,\varphi_2\}$ chart, so that \rife{Alex.intro} holds therein. In particular,
if $E$ is any open and relatively compact
Borel set in $\om_2$, then the equality is always strict, since $K$  is never constant in $E$.


\bigskip


\bigskip
\bigskip



\begin{thebibliography}{99}

\bibitem{Ambrosio}
L. Ambrosio, J. Bertrand, {\em On the regularity of Alexandrov surfaces with curvature bounded below}.
Preprint (2014).

%\bibitem{Adams} R. A. Adams, "Sobolev Spaces", Academic Press, New-York San Francisco London, 1975.

\bibitem{Ale} A. D. Alexandrov, {\em Die innere Geometrie der konvexen Fl\"{a}chen}, Springer Verlag,
Berlin, 1955.

\bibitem{Ale2} A. D. Alexandrov, V.A. Zalgaller, {\rm Intrinsic Geometry of Surfaces}, AMS Transl. Math.
Monographs, Vol. 15, Providence, RI, 1967.

\bibitem{Band0} C. Bandle, {\em On a differential Inequality and its applications
to Geometry}, {Math. Zeit.} {\bf 147}, (1976) 253-261.

\bibitem{Band} C. Bandle, {\rm Isoperimetric Inequalities and Applications}, Pitman, Boston, 1980.

%\bibitem{B1} D. Bartolucci, {\em A "Sup + C Inf" inequality for the
%equation $-\Delta u=\frac{V}{|x|^{2\alpha}} e^{u}$},
%{Proc. Royal Soc. of Edinburgh} {\bf 140A} (2010), 1119-1139.

%\bibitem{B3} D. Bartolucci, {\em A "Sup + Inf" inequality for Liouville
%type equations with weights }, Jour.  d'Analyse
%Mathematique {\bf 117} (2012), 29-46; DOI:10.1007/s11854-012-0013-7.

\bibitem{BaDo} J.L. Barbosa, M. do Carmo, {\em A Proof of a General Isoperimetric Inequality for Surfaces},
Math. Z. {\bf 162} (1978) 245-261.

\bibitem{barjga} D. Bartolucci, {\em On the best pinching constant of conformal metrics on $\mathbb{S}^2$ with
one and two conical singularities}, Jour. Geom. Analysis {\bf 23} (2013) 855-877.

\bibitem{BCast1} D. Bartolucci, D. Castorina,
{\em Self gravitating cosmic strings and the Alexandrov's inequality for Liouville-type equations}, Comm. Cont. Math. {\bf 18}(4) (2016),
1550068 (26p.).

%\bibitem{bl} D. Bartolucci, C.S. Lin, {\em Uniqueness Results for Mean Field Equations with Singular Data},
%Comm. in P. D. E. {\bf 34}(7) (2009), 676-702.

%\bibitem{BLin2} D. Bartolucci, C.S. Lin, {\em Sharp existence results for mean field equations with singular data},
%Jour. Diff. Eq. 252(7) (2012), pp. 4115-4137.

\bibitem{BLin3} D. Bartolucci, C.S. Lin, {\em Existence and uniqueness for
Mean Field Equations on multiply connected domains at the critical parameter},
{Math. Ann.}, {\bf 359} (2014), 1-44.

%\bibitem{BLT} D. Bartolucci, C.S. Lin, G. Tarantello, {\em Uniqueness and symmetry results for
%solutions of a mean field equation on ${\mathbb{S}}^{2}$ via a new bubbling phenomenon},
%{Comm. Pure Appl. Math.} {\bf 64}(12) (2011), 1677-1730.

%\bibitem{BM2} D. Bartolucci, E. Montefusco,
%{\sl On the Shape of Blow up Solutions to a Mean Field Equation},
%{Nonlinearity} {\bf 19}, (2006), {611-631}.

%\bibitem{BM3} D. Bartolucci, E. Montefusco, {\em Blow up analysis,
%existence and qualitative properties of solutions for the two
%dimensional Emden-Fowler equation with singular potential},
%M$^{2}$.A.S. {\bf 30}(18) (2007), 2309-2327.

%\bibitem{bt} D. Bartolucci, G. Tarantello, {\sl Liouville type equations with
%singular data and their applications to periodic multivortices for the
%electroweak theory}, Comm. Math. Phys. {\bf 229} (2002), 3-47.


\bibitem{Bol} G. Bol, {\em Isoperimetrische Ungleichungen f\"{u}r Bereiche auf Fl\"{a}chen},
Jber. Deutsch. Math. Verein.  {\bf 51}, (1941), 219-257.


\bibitem{bm} H. Brezis \& F. Merle,
{\sl Uniform estimates and blow-up behaviour for solutions of $-\Delta u = V(x)e^{u}$ in two dimensions},
Comm. in P.D.E.,  {\bf 16}(8,9) (1991), 1223-1253.

\bibitem{bz} J.E. Brothers \& W.P. Ziemer,
{\sl Minimal rearrangement of Sobolev functions},
J. Reine Angew. Math. {\bf 384} (1988), 153-179.

%\bibitem{BurKos} Y.D. Burago, N. N. Kosovskiĭ, {\sl The trace of BV-functions on an irregular subset},  St. Petersburg Math. J. {\bf 22} (2011), 251–-266.


\bibitem{Bur} Y.D. Burago, V.A. Zalgaller, "Geometric inequalities", Springer Ser. Sov. Math.,
Springer-Verlag Berlin Heidelberg 1988.

%\bibitem{cl} C. C. Chen \& C. S. Lin, {\sl A sharp sup+inf inequality for a nonlinear elliptic
%equation in $\rdue$,} Comm. An. Geom.,  {\bf 6}(1) (1998), 1-19.

%\bibitem{ccl} S.Y.A. Chang, C.C. Chen \& C.S. Lin, {\sl Extremal functions for a mean field equation in two dimension},
%in: "Lecture on Partial Differential Equations", New Stud. Adv. Math. {\bf 2} Int. Press, Somerville, MA, 2003, 61-93.

%\bibitem{CLin3} C.C. Chen, C.S. Lin, {\em On the Symmetry of Blowup Solutions to a
%Mean Field Equation}, Ann. Inst. H. Poincar\'e Anal. Non Lin\'eaire {\bf 18}(3) (2001), 271-296.

\bibitem{dep} De Pascale L. {\sl The Morse-Sard theorem in Sobolev spaces}, Indiana Univ. Math. J.,
{\bf 50} (2001), 1371-1386.

%\bibitem{EG} L.C. Evans, R.F. Gariepy, "Measure theory and fine properties of functions",
%CRC Press (1992).

%\bibitem{fia} F. Fiala, {\sl Le probl\`{e}me des isoperim\`{e}tres sur le surfaces ouvert \`{a} courbure
%positive}, Comm. Math. Helv., {\bf 9} (1921), 154-160.

\bibitem{Fang} H. Fang, M. Lai, {\em On curvature pinching of conic 2-spheres}, Preprint (2015). 


\bibitem{Gilb} D. Gilbarg, N.S. Trudinger, {\sl Elliptic Partial Differential Equations of
Second Order}, II Ed., Springer 2001.

%\bibitem{HK} W. K. Hayman, P.B. Kennedy {\sl Subharmonic functions}, (1076) Academic Press London-N.Y.-San Francisco

\bibitem{Hub0} A. Huber,  {\sl On the isoperimetric inequality on surfaces of variable Gaussian
curvature}, Ann. Math., {\bf 60}(2) (1954), 237-247.

\bibitem{Hub} A. Huber,  {\sl Zur Isoperimetrischen Ungleichung Auf Gekr\"{u}mmten Fl\"{a}chen},
Acta. Math.,  {\bf 97} (1957), 95-101.

%\bibitem{ls} Y.Y. Li \& I.Shafrir, {\sl Blow-up analysis for Solutions of $-\Delta u = V(x)e^{u}$
%in dimension two}, {Ind. Univ. Math. J.},  {\bf 43}(4) (1994), 1255--1270.

%\bibitem{km}
%T. Kuusi, G. Mingione, {\sl universal potential estimates}, Jour. Funct. Anal {\bf 262} (2012), 4205--4269.

%\bibitem{Lio}
%J. Liouville,
%"{\it Sur L' \'Equation aux Diff\'erence Partielles
%$\frac{d^{2} \log{\lm}}{du dv} \pm \frac{\lm}{2 a^{2}}=0$}",\\
%{J. Math. Pure Appl.} {\bf 36} 71-72 (1853).

%\bibitem{St4} W. Littman, G. Stampacchia \& H. F. Weinberger,
%{\sl Regular points for elliptic equations with discontinuous coefficients}
%Ann. Scuola Normale Sup. Pisa,  {\bf 17} (1963), 43-77.

%\bibitem{Oss1} R. Ossermann, {\em The isoperimetric inequality}, Bull. A.M.S. {\bf 84}(6) (1978), 1182-1238.

\bibitem{Oss2} R. Ossermann, {\em Bonnesen-style Isoperimetric Inequalities}, Am. Math. Mont. {\bf 86}(1) (1979), 1-29.

%\bibitem{Pic} E. Picard, {\sl De l'int\'{e}gration de l'\'{e}quation $\Delta u= e^u$ sur une surface de Riemann ferm\'{e}e},
%J. Crelle, {\bf 130}, (1905).

\bibitem{pom} Ch. Pommerenke, {\sl Boundary Behaviour of Conformal Maps}, Grandlehren der Math.
Wissenschaften, {\bf 299}, p. 300, Springer-Verlag, Berlin-Heidelberg, 1992


%\bibitem{Riesz} F. Riesz,
%{\em Sur les fonctions subharmoniques et leur rapport it la theorie du potentiel I},
%Acta Math. {\bf 48} (1926), 329-343; {\em Sur les fonctions subharmoniques et leur rapport it la theorie du potentiel II }Acta Math. {\bf 54} (1930), 321-360.

%\bibitem{Res1} Reshetnyak, Y.G. {\em Isothermal coordinates in manifolds of bounded curvature}, Doklady
%Akad. Nauk SSSR (N.S.) 94, (1954). 631-633. (Russian)

\bibitem{Res3} Y.G. Reshetnyak, {\em Isothermal coordinates in manifolds of bounded curvature II}, Sib. Math. J., vol. 1, (1960), p. 248-276.
(Russian)

\bibitem{Res2} Y.G. Reshetnyak,  {\em Two-dimensional manifolds of bounded curvature},  pp. 3-163 In: Y.G.
Reshetnyak (Ed.), Geometry IV, Encyclopaedia of Math. Sci., Vol. 70, Springer (1993).

%\bibitem{pt} J. Prajapat \& G. Tarantello, {\sl On a class of elliptic problems in $\rdue$:
%symmetry and uniqueness results}, Proc. Roy. Soc. Edinburgh, Sect. A {\bf 131} (2001), 967--985.

%\bibitem{PoT}
%A. Poliakovsky \& G. Tarantello, {\em On singular Liouville systems},
%Preprint (2013).

%\bibitem{S} I.Shafrir, {\sl A Sup+CInf inequality for the equation $-\Delta u = V(x)e^{u}$},
%{C. R. Acad. Sci. Paris},  {\bf 315}(Ser.2) (1992), 159--164.

\bibitem{suz} T. Suzuki, {\em Global analysis for a two-dimensional elliptic eiqenvalue problem with the exponential 
nonlinearly}, Ann. Inst. H. Poincar\'e Anal. Non Lin\'eaire {\bf 9}(4) (1992), 367-398.


%\bibitem{Tar} G. Tarantello, in preparation


%\bibitem{T1} G. Tarantello,
%{\sl A Harnack inequality for Liouville type equations with singular sources},
%Indiana. Univ. Math. Jour., {\bf 54}(2) (2005), 599--615.

%\bibitem{T3} G. Tarantello,
%{\sl Analytical aspects of Liouville type equations with singular sources}, Handbook Diff. Eqs., North Holland,
%Amsterdam, Stationary partial differential equations, {\bf I} (2004), 491--592.

\bibitem{stam} G.  Stampacchia, {\em Le probl{\`e}me de Dirichlet pour les {\'e}quations elliptiques
du second ordre {\`a} coefficients discontinus}, Ann. Ins. Fourier {\bf 15}(1) (1965), pp. 189-257.


\bibitem{Tal} G. Talenti {\em Best Constant in Sobolev Inequality}, Ann. Mat. Pura Appl., {\bf 110}
(1976), pp. 353-372.


\bibitem{topp} P. Topping {\em Mean curvature flow and geometric inequalities}, J. Reine Angew. Math.,
{\bf 503}, (1998) 47-61.

\bibitem{topp1} P. Topping {\em The isoperimetric inequality on a surface}, Man. Math.,
{\bf 100}, (1999) 23-33.

\bibitem{Troy0} M. Troyanov, {\em Metrics of constant curvature on a sphere with two
conical singularities}, Proc. Third Int. Symp. on Diff. Geom. (Peniscola 1988),
Lect. Notes in Math. {\bf 1410} Springer-Verlag, 296--308.

\bibitem{Troy} M. Troyanov, {\sl Prescribing curvature on compact surfaces with
conical singularities}, Trans. Amer. Math. Soc. {\bf 324} (1991), 793--821.

\bibitem{Troy1} M. Troyanov, {\em Les surfaces a courbure integrale bornee au sense d' Alexandrov}, preprint arXiv:0906.3407v1.

\bibitem{Troy2} M. Troyanov, {\em Une principe de concentration-compacit\'e pour le suites de surfaces Riemanniens},
Ann. Inst. H. Poincar\'e Anal. Non Lineaire {\bf 5} (1991), 419--441.

%\bibitem{w} G. Wolansky, {\em On steady distributions of self-attracting
%clusters under friction and fluctuations}, Arch. Rational Mech. An.
%{\bf 119} (1992), 355--391.


\end{thebibliography}
\end{document}